\documentclass[review]{elsarticle}
\pdfoutput=1

\journal{}

\usepackage{amsmath}
\usepackage{amsfonts}
\usepackage{amssymb}
\usepackage{amsthm}
\usepackage{graphicx}
\usepackage{xspace}
\usepackage{algorithm}
\usepackage{booktabs}

\usepackage[final]{listofsymbols}
\usepackage{color}
\usepackage{colortbl}
\usepackage[table]{xcolor}
\definecolor{Gray}{gray}{0.9}
\newcolumntype{g}{>{\columncolor{Gray}}c}

\usepackage{standalone}
\usepackage{pgfplots}
\usepackage{tikz}
\usetikzlibrary{arrows,automata,positioning,calc,fit,trees}

\tikzset{rectangle/.append style={draw=none,
                                  fill=gray!30,
                                  rounded corners=2pt,
                                  minimum width=1.6cm,
                                  text width=1.6cm,
                                  align=center,
                                  inner sep = 4pt,
                                  minimum height=0.4cm}}

\tikzset{rectangle_algo/.append style={draw=none,
                                  fill=gray!30,
                                  rounded corners=2pt,
                                  minimum width=2.5cm,
                                  minimum height=0.7cm,
                                  align=center,
                                  inner sep=4pt,
                                  outer sep=0pt,
                                  }}

\usepackage[list=true, labelformat=brace, position=bottom]{subcaption}
\usepackage[size=scriptsize]{todonotes}
\usepackage{multirow}
\usepackage{rotating}
\usepackage{textcomp}

\usepackage[acronym]{glossaries}
\newacronym{admm}{\textsc{admm}}{Alternating direction method of multipliers}
\newacronym{aladin}{\textsc{aladin}}{Augmented Lagrangian based Alternating Direction Inexact Newton method}
\newacronym{rapidpf}{rapid\textsc{pf}}{rapid prototyping for distributed power flow}
\newacronym{bfgs}{\textsc{bfgs}}{Broyden-Fletcher-Goldfarb-Shanno}
\newacronym{tsos}{\textsc{tso}s}{transmission system operators}
\newacronym{dsos}{\textsc{dso}s}{distribution system operators}
\newacronym{qp}{\textsc{qp}}{quadratic program}

\newcommand{\matlab}{\textsc{matlab}\xspace}
\newcommand{\papertitle}{Distributed Power Flow and Distributed Optimization -- Formulation, Solution, and\\Open Source Implementation}
\newcommand{\norm}[1]{\left\lVert#1\right\rVert}
\newcommand{\n}{n}
\newcommand{\matpower}{\textsc{matpower}\xspace}
\newcommand{\powermodels}{PowerModels\xspace}
\newcommand{\pandapower}{pandapower\xspace}
\newcommand{\aladinalpha}{\textsc{aladin}-$\alpha$\xspace}
\newcommand{\casadi}{\textsc{c}as\textsc{ad}i\xspace}
\newcommand{\fmincon}{\texttt{fmincon}\xspace}
\newcommand{\fminunc}{\texttt{fminunc}\xspace}
\newcommand{\worhp}{\texttt{worhp}\xspace}
\newcommand{\pv}{\textsc{pv}\xspace}
\newcommand{\pq}{\textsc{pq}\xspace}
\newcommand{\license}{\textsc{bsd}-3-clause\xspace}

\newtheorem{remark}{Remark}

\usepackage{enumitem}

\newlist{mylist}{enumerate}{1}
\setlist[mylist]{label*=\textsc{c}\arabic*:~, ref=\textsc{c}\arabic*}

\usepackage{hyperref}
\hypersetup{
    colorlinks=true,
    allcolors=.,
}

\begin{document}

\opensymdef
\newsym[Number of regions]{nregions}{\n^{\text{reg}}}
\newsym[Set of regions]{regions}{\mathcal{R}}
\newsym[Number of connections]{nconnections}{\n^{\text{conn}}}
\newsym[Number of core nodes]{ncore}{\n^\text{core}}
\newsym[Number of core nodes]{ncopy}{\n^\text{copy}}
\newsym[Number of buses]{nbus}{\n^\text{bus}}

\newsym[Voltage phase angle]{vang}{\phi}
\newsym[Voltage magnitude]{vmag}{v}
\newsym[Active power]{p}{p}
\newsym[Reactive power]{q}{q}

\newsym[State]{state}{x}
\newsym[State of core node]{stateCore}{x}
\newsym[State of copy node]{stateCopy}{z}

\newsym[Power flow problem]{pfproblem}{g}
\newsym[Power flow equations]{pf}{g^{\text{pf}}}
\newsym[Bus specifications]{busspecs}{g^{\text{bus}}}

\newsym[]{DistOptLocalState}{\chi}
\newsym[]{DistOptUpdatedState}{\zeta}

\closesymdef

\begin{frontmatter}

\title{\papertitle}

\author[mymainaddress]{Tillmann M\"uhlpfordt\corref{mycorrespondingauthor}\fnref{funding}}
\cortext[mycorrespondingauthor]{Corresponding author: tillmann.muehlpfordt@kit.edu}

\author[mymainaddress]{Xinliang Dai\fnref{funding}}

\author[mymainaddress]{\\Alexander Engelmann\fnref{current_address}}
\fntext[current_address]{Current address: Institute for Energy Systems, Energy Efficiency and 
Energy Economics, \textsc{tu} Dortmund, Germany}
\fntext[funding]{The authors acknowledge funding from the German Federal Ministry of Education and Research within the project \emph{MOReNet -- Modellierung, Optimierung und Regelung von Netzwerken heterogener Energiesysteme mit volatiler erneuerbarer Energieerzeugung}.}

\author[mymainaddress]{Veit Hagenmeyer}

\address[mymainaddress]{Institute for Automation and Applied Informatics\\Karlsruhe Institute of Technology\\Germany}

\begin{abstract}
  Solving the power flow problem in a distributed fashion empowers different grid operators to compute the overall grid state without having to share grid models---this is a practical problem to which industry does not have off-the-shelf answers.
  In cooperation with a German transmission system operator we propose two physically consistent problem formulations (feasibility, least-squares) amenable to two solution methods from distributed optimization (the \gls{admm}, and the \gls{aladin}); with \gls{aladin} there come convergence guarantees for the distributed power flow problem.
  In addition, we provide open source \matlab code for \gls{rapidpf}, a fully \matpower-compatible software that facilitates the laborious task of formulating power flow problems as distributed optimization problems; the code is available under \url{https://github.com/KIT-IAI/rapidPF/}.
  The approach to solving distributed power flow problems that we present is flexible, modular, consistent, and reproducible.
  Simulation results for systems ranging from 53 buses (with 3 regions) up to 4662 buses (with 5 regions) show that the least-squares formulation solved with \gls{aladin} requires just about half a dozen coordinating steps before the power flow problem is solved.
\end{abstract}


\end{frontmatter}


\section{Introduction}
The power flow problem is \emph{the} cornerstone problem for power systems analyses: find all (complex) quantities in an \textsc{ac} electrical network in steady state.
The information drawn from the solution of the power flow problem is relevant for planning power systems, as well as expanding and operating them~\cite{Grainger94book}.
Hence, the power flow problem is relevant to all stakeholders maintaining a well-functioning electrical grid, most importantly \gls{tsos} and \gls{dsos}.
Traditionally, \gls{tsos} and \gls{dsos} solve and use power flow problems independently of each other, each making modeling assumptions with respect to the other system, e.g. treating the distribution system as a lumped load for the transmission system~\cite{Sun2008}.
Clearly, these modeling assumptions---even if they were valid---may lead to real-world mismatches in both power and voltage.
Hence---what is a sovereignty-preserving way to solve power flow problems for large power systems that may be composed of several \gls{tsos} and/or \gls{dsos}?\footnote{Sometimes the literature refers to a power flow problem for a combination of \gls{tsos} and \gls{dsos} as a \emph{global power flow} problem~\cite{Sun2008, Sun2015}.\label{foot:global-power-flow}}
This is the motivating question the present paper.

To answer this question we focus on the mathematical formulation and solution of the power flow problem.
Mathematically, the power flow problem is modeled as a system of nonlinear equations, traditionally solved by Newton-Raphson methods or Gauss-Seidel approaches.
These solution techniques may be classified as \emph{centralized approaches}, i.e. the full grid model is available to a single central entity.
This entity solves the power flow problem, having access to all information not just about the problem itself but also about the solution.
Hence, this established approach is in principle able to solve power systems composed of \gls{tsos} and \gls{dsos} (so-called global power flow problems~\cite{Sun2008, Sun2015})---but only at the cost of giving up sovereignty.

Recently, so-called \emph{distributed approaches} have drawn significant academic attention.
These are methods for which several entities (or agents) solve sub-problems independently of each other, then broadcast some---but not all---information to a coordinator~\cite{Erseghe2014,Kim1997,Hug2015,Kim2000,Engelmann2019Aladin}.
The coordinator then solves a coordination problem, and sends to all entities the information they need to solve their sub-problems again.\footnote{We clearly distinguish between \emph{distributed} and \emph{decentralized} approaches, the latter requiring no central coordinator whatsoever.}
This process is repeated until convergence is achieved.
The high-level description of distributed approaches suggests several advantages relative to centralized methods:
\begin{itemize}
  \item distribute the computational effort,
  \item preserve sovereignty and/or privacy, e.g. grid models,
  \item decrease the vulnerability due to a single-point-of-failure, and
  \item add flexibility.
\end{itemize}

The interest in distributed approaches is not just academic; there exists a genuine desire by industry to leverage the advantages for real-world problems.
In Germany, for example, the horizontal connection between the four \gls{tsos}---50 Hertz, TenneT, Amprion, and TransnetBW---is based on centralized power flow solutions.
However, new legislation and the undergoing German \emph{Energiewende} toward more renewables force the German \gls{tsos} to focus on new vertical cooperation with the numerous \gls{dsos}.
For this vertical cooperation, centralized approaches are not favorable mainly due to privacy concerns: the host then combines the role of data owner and product owner, and introduces a possible single-point-of-failure.
Hence, it is \emph{not mainly} the distributed computational effort, but more the increased privacy, reliability, and flexibility that spur the interest of \gls{tsos} in distributed approaches.

Large Chinese cities are another example where the combined power flow problem for \gls{tsos} and \gls{dsos} is of relevance; in \cite{Sun2008} it is argued that many Chinese cities are operating both transmission and distribution systems, both of which are studied and operated separately however.
If a computational method were available to solve the combined power flow problem in terms of a privacy-preserving distributed problem, this would be helpful~\cite{Sun2008,Sun2015}.

In light of the above considerations the present paper contributes as follows:

\begin{mylist}
  \item \label{item:formulations} We present two mathematical formulations of distributed power flow problems as privacy-preserving and physically-consistent distributed optimization problems.
  \item \label{item:admm-aladin} We rigorously evaluate the applicability of the \acrfull{admm} and the \acrfull{aladin} to distributed power flow problems with up to several thousand buses.
  \item \label{item:rapidpf} We introduce \gls{rapidpf}: open-source \matlab code fully compatible with \matpower that allows to generate \matpower case files for distributed power flow problems tailored to distributed optimization; the code is available on \url{https://github.com/KIT-IAI/rapidPF/} under the \license license.
  \item \label{item:aladin-alpha} We extend \aladinalpha, a \matlab rapid-prototyping toolbox for distributed and decentralized non-convex optimization, to allow for user-defined sensitivities and three new solver interfaces (\fminunc, \fminunc, \worhp).
\end{mylist}
We explain our contributions relative to the state-of-the-art.

Ad \ref{item:formulations}:
The idea to solve a global power flow problem that stands for the combination of \gls{tsos} and \gls{dsos} was popularized by the works~\cite{Sun2008, Sun2015}.
Specifically,~\cite{Sun2008} coined the term ``Master-slave-splitting'' to highlight the idea that there is a master system to which several workers are connected;\footnote{We prefer the less inappropriate term ``worker'' instead of ``slave''.} also so-called ``boundary systems'' are introduced which make up the physical connection between the master and its worker~\cite{Sun2008}.
The solution of the overall power flow problem is done iteratively: initialize the boundary voltages, solve the power flow for the worker systems, then substitute the solution to the boundary system, and solve the master system.
This process is executed until the difference of voltage iterates is sufficiently small.
Any power flow solver can be used to solve the sub-problems.
Unfortunately, no convergence guarantees are provided, and the method was applied to systems with less than 200 buses.

In the follow-up work~\cite{Sun2015} a convergence analysis is carried out, but its practicability is limited due to mathematical settings---such as the implicit function theorem---that are difficult to relate to real-world criteria and/or data.
Also, the simulation results from~\cite{Sun2015} are the ones from~\cite{Sun2008}.
It hence remains unclear how well this method scales.
Unfortunately, neither \cite{Sun2008} nor \cite{Sun2015} provide plots on the actual convergence behavior of their method, or wall-clock simulation times, or the influence of different initial conditions---all of which are aspects relevant to practitioners.

In light of~\cite{Sun2008, Sun2015} the focus of the present paper is on the following:
\begin{itemize}
  \item clear distinction between problem formulation and problem solution;
  \item two different mathematical problem formulations that make no assumptions on the sub-problems (e.g. meshed grids vs. radial grids);
  \item convergence properties follow from theory of distributed optimization;
  \item reproducible numerical results for test systems with up to $\approx$ 4000 buses.
\end{itemize}

Ad \ref{item:admm-aladin}:
Recently, distributed optimization techniques have drawn attention for distributed optimal power flow problems.
It is especially \gls{admm} that finds widespread application for optimal power flow~\cite{Erseghe2014,Guo2017,Kim2000}.
However, \gls{admm} being a first-order optimization methods often converges relatively quickly to the vicinity of the optimal solution, but then takes numerous iterations to approach satisfying numerical accuracy~\cite{Boyd2011}.
In addition, \gls{admm} is known to be rather sensitive to both tuning and the choice of initial conditions~\cite{Sun2013}; line flow limits pose a significant obstacle for \gls{admm}~\cite{Erseghe2014}.
Furthermore, convergence guarantees for \gls{admm} apply to convex optimization problems, but optimal power flow is known to be non-convex.

There exist distributed optimization methods that are devised for non-convex problems, for instance \gls{aladin}.
With \gls{aladin} being a second-order method, it has access to curvature information that speed up convergence, at the expense of having to share more information among the sub-problems.
The proof-of-concept applicability of \gls{aladin} to distributed optimal power flow problems has been demonstrated in several recent works~\cite{Engelmann2017Aladin, Engelmann18b, Engelmann2019Aladin}; how to reduce the information exchange among the sub-problems is discussed in~\cite{Engelmann2020BiLevelAladin}.
Compared to \gls{admm}, \gls{aladin} has more favorable convergence properties: within a few dozen iterations, convergence to the optimal solution is achieved with satisfying numerical accuracy~\cite{Engelmann2019Aladin}.
Just like with \gls{admm}, however, tuning remains a challenge with \gls{aladin}.
Also, the largest test case to which \gls{aladin} was successfully applied is the 300-bus test case.

To summarize: both \gls{admm} and \gls{aladin} have demonstrated their potential for solving distributed optimal power flow problems.
It is fair to ask how both methods apply to distributed power flow problems---a question that has not been tackled before to the best of the authors' knowledge.
Our findings suggest that for the distributed power flow problem \gls{aladin} outperforms \gls{admm} far more significantly than it does for the optimal power flow problem (in terms of scalability, speed, performance, and tuning).
The main advantage of applying established techniques from distributed optimization to distributed power flow is that the convergence guarantees can be leveraged.

Ad \ref{item:rapidpf}:
For academic power system analyses \matpower is a mature, well-established, and widely adopted open source collection of \matlab code~\cite{Zimmerman2011matpower}.
It is not just the many computational facets that \matpower provides that make it popular (power flow and several relaxations, optimal power flow, unit commitment, etc.), but also the so-called \matpower case file format has inspired other open source packages, for instance {PowerModels.jl}~\cite{Coffrin2018PowerModels} or {PyPSA}~\cite{PyPSA}.
The \matpower case file format describes a power system with respect to its bus data, generator data, and branch data.
Additionally, there is a base MVA value for per-unit conversions, and for optimal power flow problems there is an entry on generator costs.
Based on both the popularity and the maturity of \matpower we provide glue code that solves the following laborious task: given several \matpower case files, and given connection information for these case files, construct a \matlab struct that corresponds to the mathematical problem formulation, and that is amenable to distributed optimization methods.
This glue code is called \gls{rapidpf}, and it is publicly available with a rich documentation---and full \matpower compatibility.
In addition, \gls{rapidpf} decreases the time-from-idea-to-result, it computes relevant sensitivities (gradients, Jacobians, Hessians), and it comes with post-processing functionalities.

The idea of \gls{rapidpf} is inspired by the \matlab packages TDNetGen~\cite{Pilatte2019TDNetGen} and AutoSynGrid~\cite{Sadeghian2020SynGrid}; the code for \gls{rapidpf} is hosted under \url{https://github.com/KIT-IAI/rapidPF/}.
From a first glance, TDNetGen seems to provide functionality similar to \gls{rapidpf}.
As written in the abstract, TDNetGen is \matlab code ``able to generate synthetic, combined transmission and distribution network models''~\cite{Pilatte2019TDNetGen}.
Unfortunately, TDNetGen is not as flexible as desired: there is currently no straightfoward way to generate TDNetGens so-called templates from arbitrary \matpower case files.
Also, the focus of TDNetGen is on \emph{generating} large test systems, not on \emph{solving} them.
In turn, \gls{rapidpf} allows to both generate test systems and \emph{prepare} them for solution by distributed optimization methods such as \gls{admm} and \gls{aladin}.
This preparation step must not be underestimated, because providing for an interface to distributed solvers is key in making distributed techniques more popular.

The focus of AutoSynGrid is on generating numerous test systems with similar statistical properties~\cite{Sadeghian2020SynGrid}.
Hence, AutoSynGrid is not directly comparable to either TDNetGen or \gls{rapidpf}.

Ad \ref{item:aladin-alpha}:
The recently published \matlab toolbox \aladinalpha provides several implementations of both \gls{admm} and \gls{aladin}~\cite{Engelmann2020Aladin}.
Its user interface allows the user to provide merely the cost functions, the equality constraints, and the inequality constraints.
Besides setting several default parameter settings, \aladinalpha computes derivatives required for either \gls{admm} or \gls{aladin}.
To do so, \aladinalpha relies internally on \casadi, an automatic differentiation framework that also parses the optimization problem to the low level interface of Ipopt~\cite{Andersson2019, Waechter06}.
The idea of \aladinalpha is to provide rapid prototyping capabilities for general distributed optimization problems; it is not specifically tailored to distributed (optimal) power flow problems.
From the authors' experience, this all-purpose character in combination with \casadi being hard-wired into \aladinalpha hinders it from being applicable to mid- to large-scale power flow problems.

We forked the code and tailored it to the needs of distributed power flow problems: the exact power flow Jacobian is passed from \matpower, Hessian approximations are provided, and three new solvers are interfaced (\fminunc, \fminunc, \worhp).
Also, the interface of \aladinalpha needed substantial changes to allow for passing of user-supplied sensitivities instead of auto-computed sensitivities from \casadi.

Although the motivation for the present paper is to solve power flow problems for systems composed of \gls{tsos} and \gls{dsos}, the authors stress that this setup is not a requirement.
The presented methodology is generic in the following sense:

\begin{quote}
  Given $i \in \{1, \hdots, \nregions \}$ power flow problems, and given suitable connection information, what is a coherent methodology for solving the overall power flow problem in a distributed manner?
\end{quote}

It may be that the individual power flow problems happen to coincide with \gls{tsos} and/or \gls{dsos}, but they can as well be sub-problems of a genuinely large power flow problem that should be solved in a distributed way.
In either case, the answer the present paper can provide to the above question is:

\begin{quote}
  If the distributed power flow problem is formulated as a distributed-least squares problem, and if this problem is solved with \gls{aladin} using a Gauss-Newton Hessian approximation, then the solution is found within half a dozen \gls{aladin} iterations for systems ranging from 53 to 4662 buses.
\end{quote}

\begin{remark}[Partitioning]
  The present paper assumes that the partitioning of the grid is \emph{given}.
  For insights on how to partition large grids in computationally advantageous ways we refer to~\cite{Guo2016Partitioning,Murray2019Partitioning,Kyesswa2020Partitiong}.
\end{remark}

The paper is organized as its title suggests: formulation, solution, implementation, followed by an extensive section on results, and concluding comments.
The formulation \autoref{sec:formulation} introduces nomenclature and the mathematical formulation of the distributed power flow problem.
The solution \autoref{sec:solution} covers two methods from distributed optimization: \gls{admm} and \gls{aladin}.
The implementation is covered in \autoref{sec:implementation}, with a strong focus on the open source \matlab code \gls{rapidpf}.
The results \autoref{sec:results} gives both qualitative and quantitative assessments of the approach, clearly demonstrating that the least-squares formulation in combination with \gls{aladin} is the most suitable solution approach.
Concluding comments in \autoref{sec:conclusion} close the paper.

\section{Problem formulation}
\label{sec:formulation}

Given a single-phase equivalent of a connected AC electrical network in steady state with $\nbus \in \mathbb{N}$ buses, solving the power flow problems means to solve a set of nonlinear equations such that the complex voltage and apparent power of all buses of the network is found.
The standard way to solve power flow problems is to apply a \emph{centralized} method: a single machine determines the solution, for instance, via Gauss-Seidel or Newton-Raphson methods.
An alternative is to distribute the computational effort to several machines, leading to so-called \emph{distributed} approaches.
Distributed approaches are promising because they eliminate single-point-of-failures, they better preserve privacy, their technical scale-up is easier, and they foster cooperation between transmission and distribution system operators.
The idea of \emph{distributed} power flow is to solve \emph{local} power flow problems within each subsystem, independently of each other, and to find consensus on the physical values of the exchanged power between the subsystems, see \autoref{fig:example-decomposition}.

\subsection{Nomenclature}
\label{sec:nomenclature}
\begin{figure*}
  \centering
  \begin{subfigure}{0.45\textwidth}
    \includegraphics[scale=0.45]{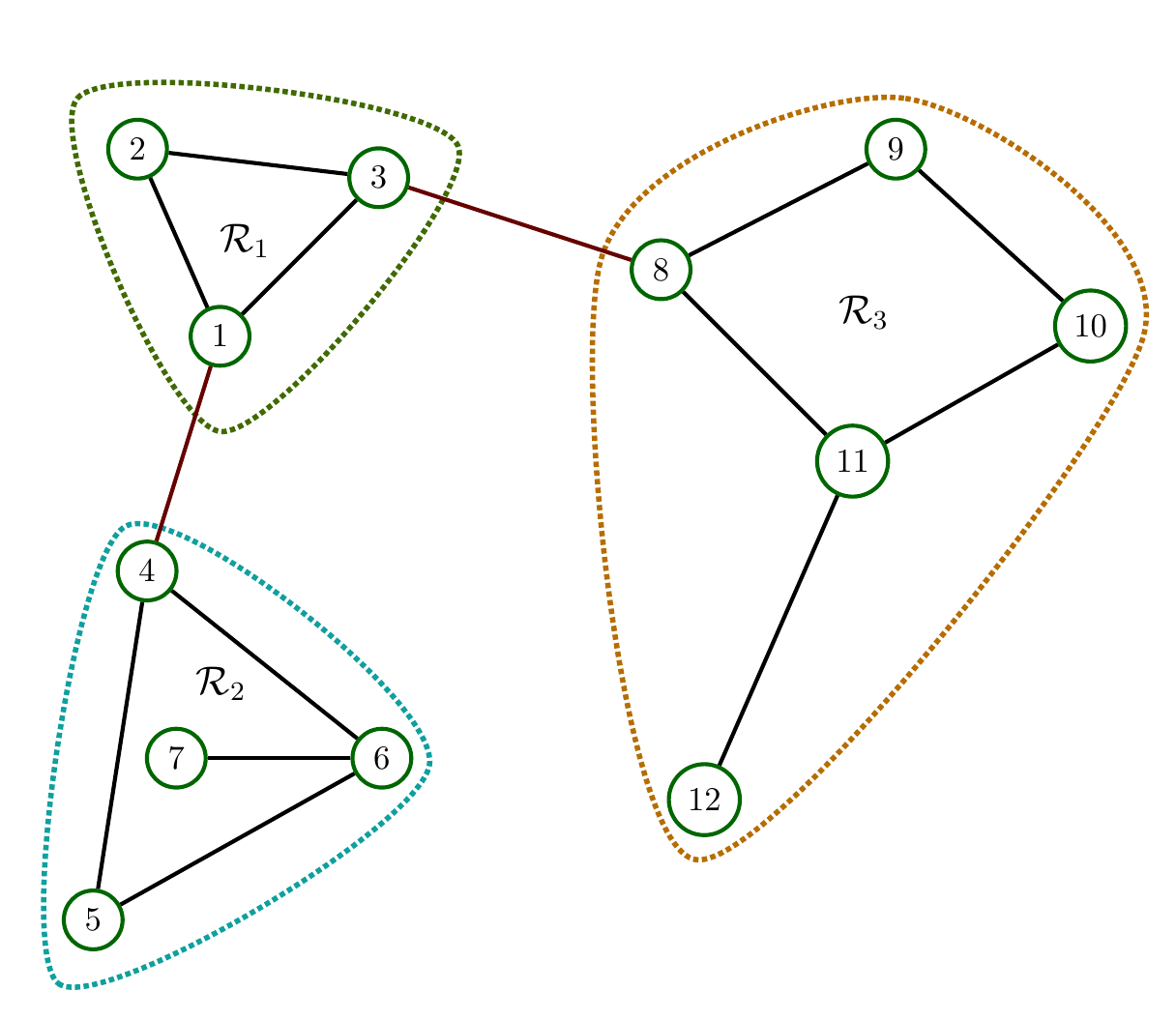}
    \subcaption{Example how to decompose a power grid into three regions $\{1, 2, 3\}$, $\{4, 5, 6, 7 \}$, and $\{8, 9, 10, 11, 12 \}$.%
    \label{fig:example-decomposition}}
  \end{subfigure}
  \hspace{0.4cm}
  \begin{subfigure}{0.45\textwidth}
    \includegraphics[scale=0.45]{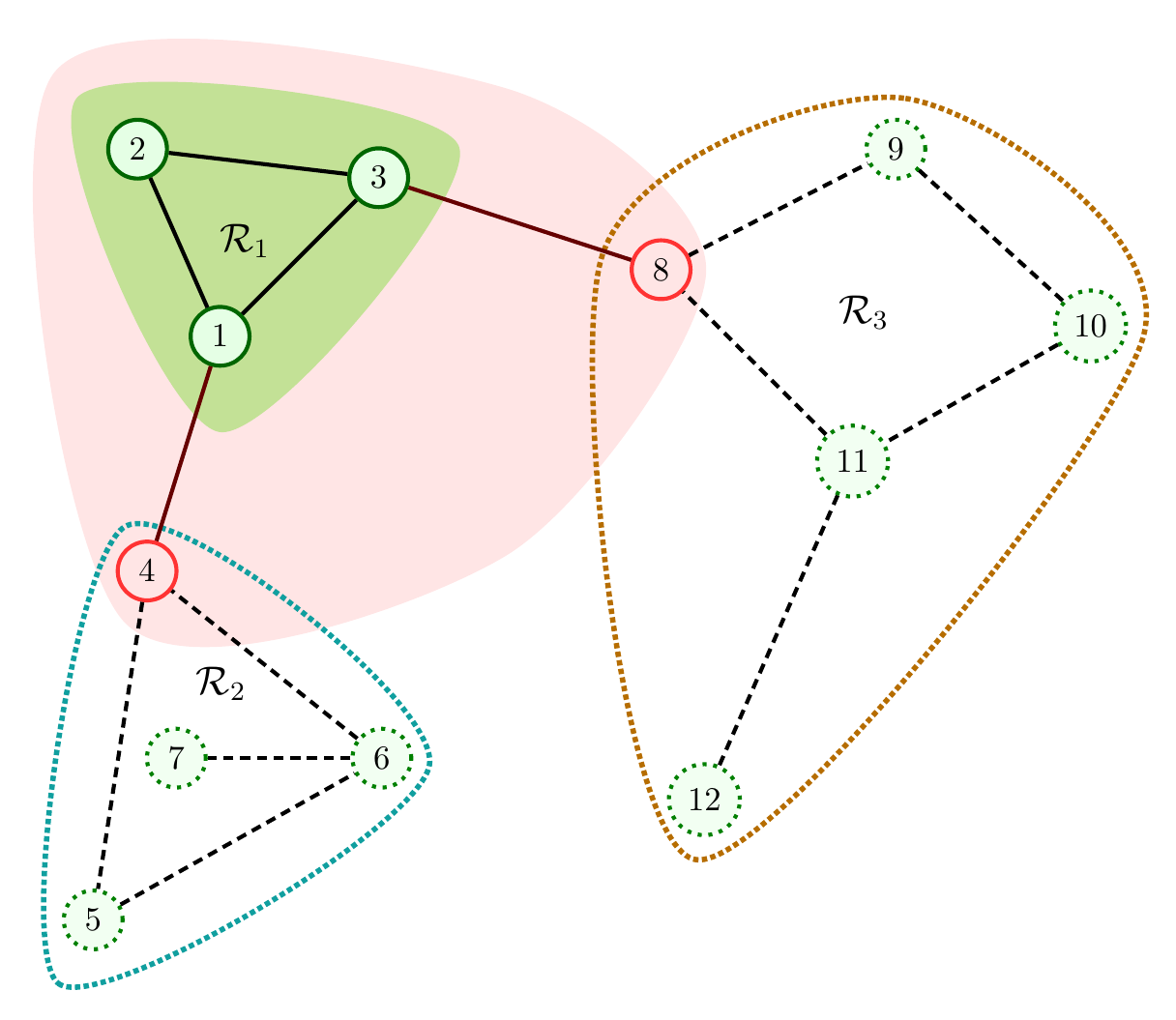}
    \subcaption{From the perspective of region $\mathcal{R}_1$, the core buses are buses $\{1, 2, 3 \}$, and the copy buses are buses $\{4, 8\}$.%
    \label{fig:example-copy-and-core-buses}}
  \end{subfigure}
  \caption{Graphical depiction of nomenclature for distributed power flow problems, see \autoref{sec:nomenclature}.}
  \label{}
\end{figure*}

Before we introduce suitable mathematical formulations for distributed power flow, we introduce some nomenclature.
For that consider \autoref{fig:example-decomposition}, which shows a 12-bus system divided into three subsystems (or so-called \emph{regions}).
Suppose we are the operator of region $\mathcal{R}_1 = \{1, 2, 3 \}$, for which we know all electrical parameters as well as all bus specifications, and for which we would like to solve a power flow problem.
This requires additional information: the complex voltages of buses $\{4, 8\}$, and the branch parameters of the tie lines---hence, connection information about the neighboring subsystems $\mathcal{R}_2 = \{4, \hdots, 7  \}$ and $\mathcal{R}_3 = \{8, \hdots, 12 \}$.\footnote{We stress that no information about the \emph{net power} of the neighboring buses $\{4, 8 \}$ is required to formulate the power flow equations.}
We shall call buses $\{ 1, 2, 3\}$ the \emph{core buses} of region~$\mathcal{R}_1$, and buses $\{4, 8\}$ the \emph{copy buses} of region~$\mathcal{R}_1$; \autoref{fig:example-copy-and-core-buses} highlights the distinction.
The combination of core buses \emph{and} copy buses allows to formulate a self-contained power flow problem for every region.

\begin{table}
  \centering
  \scriptsize
  \caption{List of symbols for distributed power flow.}
  \label{tab:list-of-symbols}
  \begin{tabular}{ll}
    \toprule
    Symbol & Meaning \\
    \midrule
    $\nregions$ & Number of regions \\
    $\nconnections$ & Number of connecting lines between regions \\
    $\ncore_{i}$ & Number of core buses in region $i$ \\
    $\ncopy_{i}$ & Number of copy buses in region $i$ \\
    \midrule
    $\stateCore_{i}$ & Electrical state of core buses in region $i$ \\
    $\stateCopy_{i}$ & Electrical state of copy buses in region $i$\\
    \bottomrule
\end{tabular}  
\end{table}

\subsection{Distributed power flow}
\autoref{tab:list-of-symbols} introduces the notation we use from here on: we consider a finite number $i \in \{1, \hdots, \nregions \}$ of regions.
The (electrical) state~$\stateCore_i$ of region~$i$ contains the voltage angles, the voltage magnitudes, the net active power, and the net reactive power of all \emph{core} buses
\begin{align}
    \label{eq:state-core}
    \stateCore_i = \begin{pmatrix}
        {\theta_i^\text{core}} &
        {v_i^\text{core}} &
        {p_i^\text{core}} &
        {q_i^\text{core}} 
    \end{pmatrix}
    \in \mathbb{R}^{4 \ncore_i}.
\end{align}
The (electrical) state~$\stateCopy_i$ of region~$i$ contains the voltage angles and the voltage magnitudes of all \emph{copy} buses
\begin{align}
    \label{eq:state-copy}
    \stateCopy_i = \begin{pmatrix}
        \theta_i^\text{copy} & 
        v_i^\text{copy}
    \end{pmatrix}
    \in \mathbb{R}^{2 \ncopy_i}.
\end{align}
Hence, each region~$i$ is represented by a total of $4 \ncore_i + 2 \ncopy_i$ real numbers.
For all core buses of region~$i$ the respective $2 \ncore_i$ power flow equations $\pf_i \colon \mathbb{R}^{4 \ncore_i} \times \mathbb{R}^{2 \ncopy_i} \rightarrow \mathbb{R}^{2 \ncore_i}$, and the respective $2 \ncore_i$ bus specifications $\busspecs_i \colon \mathbb{R}^{4 \ncore_i} \rightarrow \mathbb{R}^{2 \ncore_i}$ make up the power flow problem for this very region~\cite{Frank2016}.\footnote{The copy buses are required solely to formulate the power flow equations.}
Subtracting the number of equations from the number of decision variables gives a total of
\begin{align}
    \label{eq:dofs}
    \underbrace{4 \ncore_i + 2 \ncopy_i}_{\text{Decision vars.}} - \underbrace{2 \ncore_i}_{\text{Power flow eqns.}} - \underbrace{2 \ncore_i}_{\text{Bus specs.}} = 2 \ncopy_i
\end{align}
missing equations per region~$i$.\footnote{Note that $\textstyle\sum_{i = 1}^{\nregions} \ncopy_i \equiv 2 \nconnections$---we introduce two copy nodes for every line connecting two regions, yielding a total of $4 \nconnections$ missing equations.\label{foot:missing-equations}}

It remains to formalize the information that every \emph{copy} bus from region~$i$ corresponds to a \emph{core} bus from a neighboring region~$j \neq i$.
An example: in \autoref{fig:example-copy-and-core-buses}, bus~$4$ is a copy bus of region~$\mathcal{R}_1$, and it is a core bus of region~$\mathcal{R}_2$.
Hence, their complex voltage must be identical.

Having introduced the nomenclature we formulate distributed power flow mathematically as follows
\begin{subequations}
  \label{eq:dist-power-flow-problem}
  \begin{align}
    \label{eq:dist-power-flow-problem-pf}
    \pf_i( \stateCore_i, \stateCopy_i ) &= 0 & \forall i \in \{1, \hdots, \nregions \}\\
    \label{eq:dist-power-flow-problem-busspecs}
    \busspecs_i ( \stateCore_i ) &= 0 & \forall i \in \{1, \hdots, \nregions \}\\
    \label{eq:dist-power-flow-problem-consensus}
    \sum_{i = 1}^{\nregions} A_i \begin{bmatrix}
        \stateCore_i \\
        \stateCopy_i
    \end{bmatrix}
    &= 0.
  \end{align}
\end{subequations}
The local power flow problem for region~$i$ is given by \eqref{eq:dist-power-flow-problem-pf} and \eqref{eq:dist-power-flow-problem-busspecs}, see \autoref{rema:power-flow-equations} and \autoref{rema:bus-specifications}; the so-called consensus matrices $A_i \in \mathbb{R}^{4 \nconnections \times (4 \ncore_i + 2 \ncopy_i)}$ enforce equality of the voltage angle and the voltage magnitude at the copy buses and their respective core buses, hence they provide the remaining $4 \nconnections = \textstyle\sum_{i = 1}^{\nregions} 2 \ncopy_i$ missing equations, see \autoref{foot:missing-equations}.

\begin{remark}[Power flow equations]
  \label{rema:power-flow-equations}
  The specific form of the regional power flow equations~$\pf_i(\cdot)$ in~\eqref{eq:dist-power-flow-problem} is arbitrary.
  Nevertheless, we chose polar coordinates for the voltage phasors when defining the electrical state in~\eqref{eq:state-core}.
  In that case, the regional power flow equations are
  \begin{subequations}
    \begin{align}
      p_j &= v_j \sum_{k = 1}^{\n_i} v_k \left( G_{jk} \cos(\delta_j - \delta_k) + B_{jk} \sin(\delta_j - \delta_k) \right) \\
      q_j &= v_j \sum_{k = 1}^{\n_i} v_k \left( G_{jk} \sin(\delta_j - \delta_k) - B_{jk} \cos(\delta_j - \delta_k) \right) ,
    \end{align}
  \end{subequations}
  for all buses~$j$ from region~$i$; the bus admittance matrix entries are~$Y_{jk} = G_{jk} + \mathrm{j} B_{jk}$.
  For further details we refer to the excellent primer~\cite{Frank2012Primer}.
\end{remark}

\begin{remark}[Bus specifications]
  \label{rema:bus-specifications}
  For conventional power flow studies, each bus is modelled as one of the following:
  \begin{itemize}
    \item \emph{Slack bus:} The voltage magnitude and the voltage angle are fixed; the net active and the reactive power are determined by the power flow solution.
    \item \emph{\pq/load bus:} The active power and the reactive power are fixed; the voltage magnitude and the voltage phasor are determined by the power flow solution.
    \item \emph{\pv/voltage-controlled bus:} The active power and the voltage magnitude are fixed; the reactive power and the voltage angle are determined by the power flow solution.
  \end{itemize}
  Mathematically, these requirements are simple equality constrains of the form~$\busspecs_i(\cdot)$ for every region~$i$.
\end{remark}

\begin{remark}[Physical consistence]
  \label{rema:consistence}
  The concept of \emph{core buses} and \emph{copy buses} allows to compose the distributed power flow problem in a physically consistent manner: no additional modeling assumptions are introduced or required.
  If the solution to the distributed problem is found, then this is also the solution to the respective centralized power flow problem.
  In other words, copying buses \emph{does not} introduce a structural numerical error~\cite{Erseghe2014}.

  Other approaches, such as ``cutting'' connecting tie lines and enforcing equality of the electrical state at the intersection~\cite{Engelmann2020Aladin}, are in general \emph{not} physically consistent (only in the absence of line capacitance).
  Hence, even if the true solution to the distributed problem is found, this solution is not numerically identical to the solution of the centralized power flow problem.
  In other words, cutting lines \emph{does} introduce a structural numerical error, generally speaking.
\end{remark}

\begin{remark}[Privacy]
  \label{rema:privacy}
  To formulate the power flow equations for region~$i$, the voltage information of the copy buses needs to be shared among neighboring regions; this is inherent to the idea of \emph{core} and \emph{copy} buses.
  Although this means having to share data, the copy bus voltage data (i) does not contain a wealth of privacy information yet (ii) allows for a physically consistent problem formulation, see \autoref{rema:consistence}.
\end{remark}

\subsection{Distributed optimization problem}
The distributed power flow problem from~\eqref{eq:dist-power-flow-problem} is a system of nonlinear equations.
In contrast to the standard power flow problem, however, Problem~\eqref{eq:dist-power-flow-problem} is in a form amenable to distributed optimization.
Specifically, we propose to solve Problem~\eqref{eq:dist-power-flow-problem} either as a \emph{distributed feasibility problem}
\begin{subequations}
    \label{eq:dist-feasibility-problem}
    \begin{align}
      \underset{\underset{\forall i \in  \{1, \dots, \nregions\}}{\stateCore_i, \stateCopy_i }}{\operatorname{min}} \: 0 \quad \operatorname{s.t.}\\
        \pf_i( \stateCore_i, \stateCopy_i ) &= 0 \label{eq:pfeq}\\
        \busspecs_i ( \stateCore_i ) &= 0 \label{eq:busSpec} \\
        \sum_{i = 1}^{\nregions} A_i \begin{bmatrix}
            \stateCore_i \\
            \stateCopy_i
        \end{bmatrix}
        &= 0,
    \end{align}
\end{subequations}
or as a \emph{distributed least-squares problem}\footnote{If not stated otherwise, we have $\| \cdot \| \equiv \| \cdot \|_2$.}
  \begin{align}
    \label{eq:dist-least-squares-problem}
    \underset{\underset{\forall i \in  \{1, \dots, \nregions\}}{\stateCore_i, \stateCopy_i }}{\operatorname{min}} \sum_{i = 1}^{\nregions}  \: \norm{\begin{bmatrix}
      \pf_i( \stateCore_i, \stateCopy_i ) \\
      \busspecs_i ( \stateCore_i )
    \end{bmatrix}}^2  \quad \operatorname{s.t.} \quad
         \sum_{i = 1}^{\nregions} A_i \begin{bmatrix}
            \stateCore_i \\
            \stateCopy_i
        \end{bmatrix}
      = 0.
  \end{align}
Necessarily, the solution from the distributed feasibility problem is a solution for the distributed least-squares problem.

To summarize, we propose to formulate distributed power flow problems as either a distributed feasibility problem~\eqref{eq:dist-feasibility-problem} or a distributed least-squares problem~\eqref{eq:dist-least-squares-problem}.
Both formulations divide and conquer: formulate power flow problems for every region, and relate them by enforcing equal voltages at the connecting buses.
The privacy overhead for the regional power flow problems is limited: only the voltage information of the connecting buses is required to formulate the regional power flow equations.

Both of the given formulations---feasibility~\eqref{eq:dist-feasibility-problem} and least-squares~\eqref{eq:dist-least-squares-problem}---are special cases of a more general problem formulation.
We shall state the general problem formulation in order to simplify the solution algorithms to follow.
Using
\begin{subequations}
\begin{equation}
  \regions = \{1, \hdots, \nregions \},
\end{equation}
we define
  \label{eq:sepForm}
  \begin{align} 
    \underset{\underset{\forall i \in \regions}{\DistOptLocalState_i }}{\operatorname{min}}~
 \sum_{i \in \regions}  f_i(\DistOptLocalState_i) \\
  \text{subject to }\quad  g_i(\DistOptLocalState_i)&=0 \qquad \forall i \in  \regions
  \label{eq:sepProbGi} \\
  \sum_{i \in \regions} A_i \DistOptLocalState_i &= 0,
  \label{eq:consConstr}
	\end{align}
\end{subequations}
where $\DistOptLocalState_i = ( \stateCore_i, \stateCopy_i )$ combines the core bus state and the copy bus state for region~$i$.
The consensus constraints~\eqref{eq:consConstr} are identical for either problem formulation; the correspondence of the cost and the equality constraints is summarized in the following \autoref{tab:correspondences}.

\begin{table}[h!]
  \scriptsize
  \centering
  \caption{Correspondence of terms from general problem~\eqref{eq:sepForm} to feasibility problem~\eqref{eq:dist-feasibility-problem} and to least-squares problem~\eqref{eq:dist-least-squares-problem}.}
  \label{tab:correspondences}
  \begin{tabular}{ccc}
    \toprule
    Term from~\eqref{eq:sepForm}& Feasibility problem~\eqref{eq:dist-feasibility-problem} & Least-squares problem~\eqref{eq:dist-least-squares-problem} \\
    \midrule
    $f_i(\DistOptLocalState_i) = $ & 0 & $\norm{\begin{bmatrix}
      \pf_i( \stateCore_i, \stateCopy_i ) \\
      \busspecs_i ( \stateCore_i )
    \end{bmatrix}}^2$ \\
    $g_i(\DistOptLocalState_i) = $ & $\begin{bmatrix}
      \pf_i( \stateCore_i, \stateCopy_i ) \\
      \busspecs_i ( \stateCore_i )
    \end{bmatrix} $ & n/a \\
    \bottomrule
  \end{tabular}
\end{table}

\begin{remark}[Nonlinear least-squares problems~\cite{Nocedal2006}]
  \label{rema:nonlinear-least-squares}
  Least-squares have been studied extensively.
  Besides being an intuitive formulation of a problem at hand, least-squares problems provide rich structure that can be exploited.
  For nonlinear least-squares problems, it is well-known that Gauss-Newton methods work well.
  Instead of applying a full Gauss-Newton method we merely exploit the fact that the Hessian matrix can be approximated by the matrix product of the Jacobian.
\end{remark}

\section{Problem solution}
\label{sec:solution}

Two viable methods to tackle distributed optimization problems of the form~\eqref{eq:dist-feasibility-problem} or~\eqref{eq:dist-least-squares-problem} are \gls{admm} and \gls{aladin}; we provide a brief  overview of both.
In the following, the superscript \textsuperscript{$k$} denotes the $k$\textsuperscript{th} iterate; the superscript \textsuperscript{$0$} hence denotes the initial condition.

\begin{remark}[Wording]
  \label{rema:wording}
  Different problems bring about different wording.
  In the problem formulation in \autoref{sec:formulation} we speak of ``regions'', because the power flow problem is usually related to an existing physical region.
  In the problem solution to follow, however, we prefer to speak of ``subsystems'', and ``local problems'', because the optimization problems that need to be solved in parallel need not resemble anything that exists in the physical world.
\end{remark}

\subsection{\gls{admm}}
\label{sec:solution-admm}
The \acrfull{admm} is a popular method for distributed optimization, particularly so for problems in context of power systems \cite{Erseghe2014,DallAnese2013,Guo2017}.
Although \gls{admm} works often well in context of power systems, convergence is in general not guaranteed due to the non-convex \textsc{ac} power flow equations.
We use \gls{admm} as a benchmark method reflecting the current state-of-the-art for distributed optimization in power systems.
Note that there exists a plethora of \gls{admm} variants; the present paper relies on the formulation from~\cite{Houska2016}.
We refer to~\cite{Houska2016,Boyd2011,Bertsekas1989} for more details on \gls{admm} and its derivations and restrict ourselves to recalling the overall algorithm.

For problem~\eqref{eq:sepForm}, \gls{admm} is summarized in Algorithm~\ref{alg:ADMMspec}.
In step~\ref{item:admm-parallel-step}, \gls{admm} solves local optimization problems, where the influence of the neighboring regions is incorporated via auxiliary terms in the objective function considering Lagrange-multiplier estimates $\lambda_i \in \mathbb{R}^{\nconnections}$  and estimates of the primal variables  $\DistOptUpdatedState_i^k \in \mathbb{R}^{4 \ncore_i + 2 \ncopy_i}$.
In step~\ref{item:admm-centralized-step}, all local solutions $\DistOptLocalState_i^{k+1}$ are collected in a coordination problem.
In many cases, this step can be simplified to a simple averaging step between neighboring subsystems~\cite{Boyd2011}.
Finally, the solution of the coordination problem $\DistOptUpdatedState_i^k$ is sent to each subsystem, and after a local Lagrange multiplier update---step~\ref{item:admm-Lagrange-update}---the iterates start from the beginning.

Advantages of \gls{admm} are its relatively small communication overhead, i.e. the necessity to communicate only the local solution between neighboring subsystems, and a relatively simple coordination, which consists of computing a simple average between neighboring subsystems.

\begin{algorithm} 
	\footnotesize
	\caption{\gls{admm} for problem \eqref{eq:sepForm}}
	\textbf{Initialization: $\DistOptUpdatedState_i^0,\lambda_i^0$ for all $i \in \mathcal{R}$, $\rho$}\\
	\textbf{Repeat:}
	\begin{enumerate}[label=\arabic*)]
		\item $\DistOptLocalState_i^{k+1}=\underset{g_i(\DistOptLocalState_i)=0}{\operatorname{argmin}} \;  f_i(\DistOptLocalState_i) + \lambda_i^{k\top} A_i \DistOptLocalState_i  + \frac{\rho }{2}\|A_i (\DistOptLocalState_i -\DistOptUpdatedState_i^k)\|_2^2, \;\;\, i \in \mathcal{R}$ \hfill (parallel)
		\label{item:admm-parallel-step}
		\item $\DistOptUpdatedState^{k+1}=\underset{A\DistOptUpdatedState=0}{\operatorname{argmin}} \; \sum_{i \in  \mathcal{R}} - \lambda_i^{k\top} A_i \DistOptUpdatedState_i  + \frac{\rho }{2}\|A_i (\DistOptLocalState_i^{k+1} -\DistOptUpdatedState_i)\|_2^2$\hfill (centralized)
		\label{item:admm-centralized-step}
		\item $\lambda_i^{k+1} = \lambda^k_i + \rho A_i(\DistOptLocalState_i^{k+1}-\DistOptUpdatedState_i^{k+1}),\qquad  \qquad\qquad\qquad\qquad\qquad i\in \mathcal{R}$ \hfill (parallel)
		\label{item:admm-Lagrange-update}
	\end{enumerate}
	\label{alg:ADMMspec}
\end{algorithm} 

\begin{remark}[Local optimization problems]
	\label{rema:local-optimization-problems}
	Note that the local optimization problems in step~\ref{item:admm-parallel-step} reflect the overall problem structure: in case of the distributed feasibility problem formulation~\eqref{eq:dist-feasibility-problem}, the local optimization problems comprise local feasibility problems, i.e. optimization problems with a zero cost and non-zero equality constraints; in case of the least-squares problem formulation~\eqref{eq:sepForm}, the local problems are unconstrained nonlinear least-squares problems.
\end{remark}

\begin{remark}[\gls{admm} for non-convex problems]
	Recently, the convergence of \gls{admm} has been shown in~\cite{Hong2016,Wang2019} for special classes of non-convex problems.
	However, these works consider non-convexities in the objective function, whereas in case of the \textsc{ac} power flow equations the non-convexity appears in the constraints, for which to the best of our knowledge no convergence guarantee exists so far.
	Note that divergence of \gls{admm} can occur also for very small-scale problems in context of power systems~\cite{Christakou2017}; however, this is rarely observed.
\end{remark}

\subsection{\gls{aladin}}
\label{sec:solution-aladin}
As an alternative to \gls{admm}, the \acrfull{aladin} has been proposed~\cite{Houska2016}.
Its main idea is to replace the relatively simple coordination step in \gls{admm} with a more sophisticated one including also constraint and curvature information to yield fast and guaranteed convergence---also for problems with non-convex constraints.

\gls{aladin} for problem~\eqref{eq:sepForm} is shown in Algorithm~\ref{alg:ALADIN}.
Step~\ref{item:aladin-parallel-step} of \gls{aladin} is similar to \gls{admm}: each subsystem minimizes its own objective function with auxiliary terms. 
A minor difference is that \gls{aladin} maintains one global Lagrange multiplier~$\lambda$ only and that positive definite weighting matrices~$\Sigma_i$ are considered in the augmentation term, where $\|x\|_\Sigma^2=x^\top \Sigma\, x$.
In step~\ref{item:aladin-computation-step}, \gls{aladin} then computes sensitivities of the local problems, i.e. the gradient of the cost function, $\nabla f_i^{\phantom{i}}(\DistOptLocalState_i^k)$,  an approximation of the Hessian of the Lagrangian function, $B_i^k$, and the Jacobian matrix of the constraints, $\nabla g_i(\DistOptLocalState_i^k)$.
These sensitivities are communicated to a central coordinator, which solves an equality-constrained \gls{qp} in step~\ref{item:aladin-centralized-step} of \gls{aladin}. 
As a result, primal increments $\Delta \DistOptLocalState_i^k$ are communicated back to the subsystems, which update $\DistOptUpdatedState_i^k$ and $\lambda^k$ in step~\ref{item:aladin-update-step}, and the iteration starts from the beginning. 

In \gls{aladin}, there are two tuning parameters: $\nu$ and $\rho$. 
The scaling matrices $\Sigma_i$ can be used for variable scaling---in case of well-behaved problems they can simply be set to the identity matrix. 
For details on selecting these parameters in context of power systems we refer to~\cite{Engelmann2019Aladin}.

\begin{algorithm}[h!]
    \footnotesize
    \caption{\gls{aladin} for problem \eqref{eq:sepForm}}
    \textbf{Initialization: $\DistOptUpdatedState_i^0,\lambda^0$, $\Sigma_i\succ 0$ for all $i \in \mathcal{R}$, $\nu,\rho$}\\
    \textbf{Repeat:}
    \begin{enumerate}[label=\arabic*)]
        \item Solve  for all $i \in \mathcal{R}$ 
        \[ \label{step:locMinStep}
        \phantom{\Delta} \DistOptLocalState_i^{k}=\underset{g_i(\DistOptLocalState_i)=0}{\operatorname{argmin}}\;   f_i(\DistOptLocalState_i) + \lambda^{k\top} A_i \DistOptLocalState_i + \frac{\nu}{2}\|\DistOptLocalState_i-\DistOptUpdatedState^k_i\|^2_{\Sigma}, \qquad  \qquad \qquad \quad   \text{(parallel)}\] 
        \label{item:aladin-parallel-step}
        \item Compute   $\nabla f_i^{\phantom{i}}(\DistOptLocalState_i^k)$, $ \;B_i^k\approx \nabla_{\DistOptLocalState_i}^2 \left ( f_i^{\phantom{i}}(\DistOptLocalState_i^k) + \gamma_i^\top g_i^{\phantom{i}}(\DistOptLocalState_i^k) \right )$,   $\; \nabla  g_i(\DistOptLocalState_i^k)$. 
        \label{item:aladin-computation-step}
        \item Solve the  coordination \gls{qp} \\
        \[
        \begin{aligned}
        \Delta \DistOptLocalState^k =&\; \underset{\Delta \DistOptLocalState}{\operatorname{argmin}}   \textstyle \sum_{i \in \mathcal{R}} \frac{1}{2} \Delta \DistOptLocalState_i^\top B_i^k \Delta \DistOptLocalState_i + \nabla f_i^\top(\DistOptLocalState_i^k) \Delta \DistOptLocalState_i \qquad \qquad  \qquad \;\; \text{(centralized)} \\
        &+\lambda^{k\top}\left( \textstyle \sum_{i \in \mathcal{R}} A_i(\DistOptLocalState_i^k + \Delta \DistOptLocalState_i) - b\right ) + \frac{\rho}{2} \|\textstyle \sum_{i \in \mathcal{R}} A_i(\DistOptLocalState_i^k + \Delta \DistOptLocalState_i) - b \|_2^2 \\
        &\text{subject to}\;\;  \nabla  g(\DistOptLocalState^k) \Delta \DistOptLocalState  = 0. \qquad \qquad  \qquad \qquad  \qquad \qquad  \quad \;\;\;    
        \end{aligned}
        \]
        \label{item:aladin-centralized-step}
        \item 
        Set 
        $
        \DistOptUpdatedState^{k+1}_i = \DistOptLocalState^k_i + \Delta \DistOptLocalState^k_i$ and $\lambda^{k+1}= \lambda^k + \rho \left ( \sum_{i \in \mathcal{R}} A_i\DistOptLocalState_i - b\right).
        $ \hfill (parallel)
        \label{item:aladin-update-step}
    \end{enumerate} \label{alg:ALADIN}
\end{algorithm}

\begin{remark}[Choosing Hessian approximations $B_i^k$]
\gls{aladin} is guaranteed to converge locally for any positive definite Hessian approximation $B_i^k$~\cite{Houska2016}.
However, the domain of local convergence and especially the convergence rate depend on the quality of the approximation.
Different Hessian approximations may also reduce the communication and computation overhead: for example, a \gls{bfgs} approach is chosen in~\cite{Engelmann2019Aladin}.
We study the influence of the choice of the Hessian approximation in~\autoref{sec:different-Hessian-approximations}.
Clearly, the choice of the Hessian should be motivated by the structure and nature of the local optimization problem, cf. \autoref{rema:local-optimization-problems}.
\end{remark}

\begin{remark}[Communication and coordination effort in \gls{aladin} and \gls{admm}]
In contrast to \gls{admm}, \gls{aladin} requires more communication and coordination per iteration compared with \gls{admm}.
The sensitivities~$B_i^k$ and~$\nabla g_i(\DistOptLocalState_i^k)$ have to be communicated, whereas \gls{admm} requires to communicate local decision variables~$\DistOptLocalState_i^k$ only.
Also the coordination step in \gls{aladin} is more expensive: instead of computing simple averages, the coordination \gls{qp} in step~3) of \gls{aladin} requires solving a linear system of equations.
With the help of this additional information, however, \gls{aladin} converges faster than \gls{admm}, hence partially compensating for the additional communication overhead.
As an alternative to basic \gls{aladin}, which is proposed here, one might consider using bi-level \gls{aladin}~\cite{Engelmann2020BiLevelAladin}.
Bi-level \gls{aladin} is able to reduce the per-step communication and coordination overhead even further.
We refer to~\cite{Engelmann2019Aladin,Engelmann2020BiLevelAladin} for more detailed analytical and numerical comparisons.
\end{remark}

\section{Implementation}
\label{sec:implementation}

The problem formulations (\autoref{sec:formulation}) and suggested solutions (\autoref{sec:solution}) are moot without means to actually implement, execute, and validate them.
We introduce \gls{rapidpf}, an open source \matlab code that tackles the problem formulation, and we present an extension to \aladinalpha, an open source \matlab code that deals with the problem solution.

\subsection{\Acrfull{rapidpf}}
Although there exist several excellent open-source tools to model, study, and solve (optimal) power flow problems (e.g. \matpower in \matlab~\cite{Zimmerman2011matpower}, \powermodels in Julia~\cite{Coffrin2018PowerModels}, or \pandapower in Python~\cite{pandapower.2018}), the same cannot be said for \emph{distributed} (optimal) power flow problems---to the best of the authors' knowledge.
To help overcome both the tedious, error-prone, and laborious task of formulating distributed power flow problems, and of interfacing distributed optimization methods, we provide open source \matlab code for \gls{rapidpf}, which automates the following task:
\begin{figure*}
  \scriptsize
  \centering
  \begin{tikzpicture}[->,>=stealth',shorten >=1pt,auto, node distance = 2mm and 5mm]
  \pgfdeclarelayer{bg} 
  \pgfsetlayers{bg, main} 
  \node[rectangle, fill=blue!20](connection_information){Connection information};
  \node[rectangle, fill=teal!20](case_file_1)[below = of connection_information]{Case file 1}; 
  \node[rectangle, fill=teal!20](case_file_2)[below = of case_file_1]{Case file 2};
  \node[rectangle, fill=teal!20](case_file_dots)[below = of case_file_2]{\dots};
  \node[rectangle, fill=teal!20](case_file_n)[below = of case_file_dots]{Case file \nregions};
  \node[rectangle, fill=gray!40](case_file_generator)[right = of case_file_2]{Case file generator};
  \node[rectangle, fill=gray!40](case_file_splitter)[right = of case_file_generator]{Case file splitter};
  \node[rectangle, fill=gray!40](case_file_parser)[right = of case_file_splitter]{Case file parser};
  \node[rectangle, fill=orange!20](aladin)[right = of case_file_parser, yshift=-9mm]{Distributed least-squares problem~\eqref{eq:dist-least-squares-problem}};
  \node[rectangle, fill=orange!20](admm)[right = of case_file_parser, yshift=9mm]{Distributed feasibility problem~\eqref{eq:dist-feasibility-problem}};

  \begin{pgfonlayer}{bg}    
    \node[draw, draw=none, rounded corners, fill=gray!10, fit={($(case_file_generator.south west)+(-2.5mm, -6mm)$) ($(case_file_parser.north east)+(2.5mm, 6mm)$)}, inner sep=0pt, label=above:Open source \matlab code \gls{rapidpf}] {};
  \end{pgfonlayer}

  \draw[->] (connection_information.east) to [](case_file_generator);
  \draw[->] (case_file_1.east) to 
  (case_file_generator);
  \draw[->] (case_file_2.east) to (case_file_generator.west);
  \draw[->] (case_file_dots.east) to (case_file_generator);
  \draw[->] (case_file_n.east) to (case_file_generator);
  \draw[->] (case_file_generator.east) to (case_file_splitter.west);
  \draw[->] (case_file_splitter.east) to (case_file_parser.west);
  \draw[->] (case_file_parser) to (admm.west);
  \draw[->] (case_file_parser) to (aladin.west);

  \end{tikzpicture}
  \caption{Flow chart for \gls{rapidpf} depicting its inputs (case files \& connection information) and its output (\matlab struct compatible with \aladinalpha).}
  \label{fig:rapidpf-flow-chart}
\end{figure*}
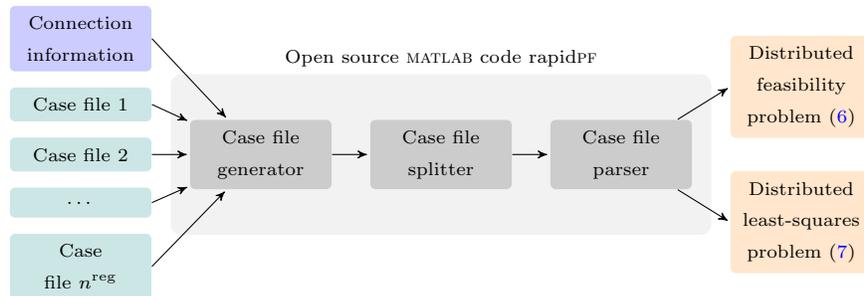
  \begin{quote}
    Given $\nregions$ \matpower case files for all $i \in \{1, \hdots, \nregions \}$ regions, and given information about how the $i \in \{1, \hdots, \nregions \}$ regions are connected, generate a \matlab struct compatible with \aladinalpha.
  \end{quote}
The features of \gls{rapidpf} span:
\begin{itemize}
  \item \emph{Rapid prototyping:} \gls{rapidpf} decreases the time-from-idea-to-result.
  \item \emph{Compatibility:} \gls{rapidpf} is compatible with \matpower and \aladinalpha. All generated case files can be visualized, for example, with the excellent ``Steady-State AC Network Visualization in the Browser''\footnote{Available on \url{https://immersive.erc.monash.edu/stac/}.}.
  \item \emph{Comparability:} \gls{rapidpf} generates \matpower case files that can be validated by \matpower functions such as \texttt{runpf()}.
  \item \emph{Sensitivities:} \gls{rapidpf} generates function handles for gradients, Jacobians, and Hessians.
  \item \emph{Documentation:} \gls{rapidpf} comes with a self-contained and user-friendly online documentation.
  \item \emph{Open source:} \gls{rapidpf} is publicly available under the \license license on \url{https://github.com/KIT-IAI/rapidPF/}.
  \item \emph{Post-processing:} \gls{rapidpf} provides rich post-processing functionalities to analyze the results quickly and intuitively.
\end{itemize}

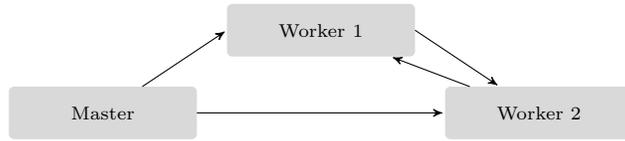
\begin{figure}
  \centering
  \scriptsize
  \begin{tikzpicture}[->,>=stealth',shorten >=1pt,auto, node distance = 4mm and 4mm, minimum width=1.0cm]
    \node[rectangle_algo](master){Master};
    \node[rectangle_algo, draw=none, fill=white](dummy)[right = of master]{};
    \node[rectangle_algo](worker_2)[right = of dummy]{Worker 2};
    \node[rectangle_algo](worker_1)[above=of dummy]{Worker 1};

    \draw[->] (master) to (worker_1.west);
    \draw[->] (master) to (worker_2);
    \draw[->] (worker_1.east) to (worker_2);
    \draw[->] (worker_2) to (worker_1);
  \end{tikzpicture}
  \caption{Supported types of connections between regions.}
  \label{fig:rapidpf-supported-connections}
\end{figure}

The code of \gls{rapidpf} is made up of three components: the case file generator, the case file splitter, and the case file parser, see \autoref{fig:rapidpf-flow-chart}.
The case file generator requires as inputs several \matpower case files in combination with their connection information; the connection information encodes \emph{who} is connected to \emph{whom} and \emph{by what} (kind of branch and/or transformer).
The regions can be connected in (almost) arbitrary ways, see \autoref{fig:rapidpf-supported-connections}.\footnote{The exception being that two buses are allowed to be connected by just one line. \autoref{rema:connecting-buses} provides further guidance about the assumptions on how buses can be connected.}
The output of the case file generator is a \matpower-compatible merged case file.
This merged case file is generated for validation purposes: it provides a reference solution that can be computed, for instance, by running \matpower's \texttt{runpf()} command.
The splitter adds information to each of the $\nregions$ case files about its core buses and copy buses.
Finally, the parser takes the augmented case files, and generates an \aladinalpha-compatible \matlab struct that describes the problem either as a distributed feasibility problem~\eqref{eq:dist-feasibility-problem} or as a distributed least-squares problem~\eqref{eq:dist-least-squares-problem}.
The parser also generates sensitivities of the power flow problem, namely the Jacobian of the power flow equations and bus specifications as well as their Hessian information.

\begin{remark}[Sensitivities]
  \label{rema:sensitivities}
  All first- and second-order optimization methods require information about derivatives.
  Hence, \gls{rapidpf} provides them for the user.
  The gradient of the local cost function, and the Jacobian of the power flow problem are the exact analytical expressions.
  The Hessian matrix---required only for \gls{aladin} but not \gls{admm}---is approximated by one of four methods: finite differences, \gls{bfgs}, limited-memory \gls{bfgs}, or Gauss-Newton.
  The first three methods can be applied to both problem formulations (feasibility~\eqref{eq:dist-feasibility-problem} and least-squares~\eqref{eq:dist-least-squares-problem}); Gauss-Newton is a method tailored to nonlinear least-squares problems~\cite{Nocedal2006}, hence applies only to the least-squares formulation~\eqref{eq:dist-least-squares-problem}.
\end{remark}

\begin{remark}[Connecting buses]
  \label{rema:connecting-buses}
  A few more words are appropriate about how systems can be connected within the case file generator.
  First, we formally distinguish between the \emph{master system} and its \emph{worker systems}.
  The sole difference is that (without loss of generality) the slack bus of the overall system is the slack bus of the master system.
  The connection between two system is directed, imposing a natural distinction between the \emph{from}- and \emph{to}-system.
  For instance, consider the line connecting the \emph{Master} and \emph{Worker 1} in \autoref{fig:rapidpf-supported-connections}: the \emph{Master} is the \emph{from}-system, \emph{Worker 1} is the \emph{to}-system.
  The connecting buses in both the \emph{from}- and the \emph{to}-system must be generation buses, hence either a slack bus or a \pv bus.
  If the connecting bus in the \emph{to}-worker-system is the slack bus, then this slack bus is replaced by a \pq bus with zero generation and zero demand.
  If the connecting bus in the \emph{to}-worker-system is a \pv bus, then this \pv bus is replaced by a \pq bus with zero generation and its original demand.
  If no connecting bus in the \emph{to}-worker-system is the slack bus, then the worker system's slack bus is replaced by a \pv bus; the respective set points for the active power and the voltage magnitude are taken from the \matpower case file entries in \texttt{mpc.gen}.
\end{remark}

\subsection{Extensions to \aladinalpha}
Whereas \gls{rapidpf} is \matlab code tailored to simplify and streamline the \emph{problem formulation}, the open source \matlab code \aladinalpha is used to tackle the \emph{problem solution}~\cite{Engelmann2020Aladin}.
\aladinalpha provides tested implementations and several variants of both \gls{admm} and \gls{aladin}.
Under the hood, \aladinalpha depends to a large degree on \casadi---an open source tool for algorithmic differentiation---and Ipopt as the solver for nonlinear programs.
Unfortunately, the sole dependency on \casadi and Ipopt hinders distributed methods from \aladinalpha to be applicable to medium- to large-scale power systems (as we shall discuss in \autoref{sec:results}).
Hence, we created a separate branch for \aladinalpha that allows to use the user-defined sensitivities from \gls{rapidpf}, and that allows to interface different solvers such as \fmincon, \fminunc,\footnote{The solvers \fmincon and \fminunc are part of \matlab's Optimization Toolbox\texttrademark.} or \worhp~\cite{Kuhlmann2018WORHP}, see also the right-hand side of \autoref{fig:formulation-solution-solvers}.

\section{Results}
\label{sec:results}
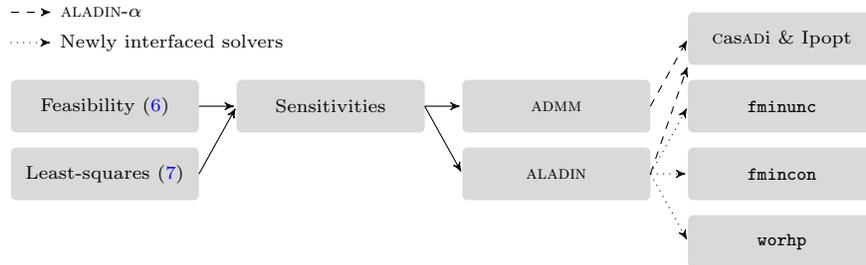
\begin{figure*}
  \scriptsize
  \centering
  \begin{tikzpicture}[->,>=stealth',shorten >=1pt,auto, node distance = 2mm and 5mm]
  \pgfdeclarelayer{bg} 
  \pgfsetlayers{bg, main} 
  \node[rectangle_algo](feasibility){Feasibility~\eqref{eq:dist-feasibility-problem}};
  \node[rectangle_algo](least-squares)[below = of feasibility]{Least-squares~\eqref{eq:dist-least-squares-problem}};

  \node[rectangle_algo](sensitivities)[right = of feasibility]{Sensitivities};

  \node[rectangle_algo](admm)[right = of sensitivities]{\gls{admm}};
  \node[rectangle_algo](aladin)[below = of admm]{\gls{aladin}};
  \node[rectangle_algo](fminunc)[right = of admm]{\fminunc};
  \node[rectangle_algo](casadi)[above = of fminunc]{\casadi \& Ipopt};
  \node[rectangle_algo](fmincon)[below = of fminunc]{\fmincon};
  \node[rectangle_algo](worhp)[below = of fmincon]{\worhp};

\draw[->] (feasibility.east) to (sensitivities.west);
\draw[->] (least-squares.east) to (sensitivities.west);
  \draw[->] (sensitivities.east) to (admm.west);
  \draw[->] (sensitivities.east) to (aladin.west);

  \draw[->, dashed] (admm.east) to (casadi.west);
  \draw[->, dashed] (aladin.east) to (casadi.south west);  
  \draw[->, dotted] (aladin.east) to (fminunc.west);
  \draw[->, dotted] (aladin.east) to (fmincon.west);
  \draw[->, dotted] (aladin.east) to (worhp.west);

  \draw[->, dashed] (current bounding box.north west)++(0,  0em) -- ++(2em, 0) node[right] {\aladinalpha};
  \draw[->, dotted] (current bounding box.north west)++(0, -2em) -- ++(2em, 0) node[right] {Newly interfaced solvers};

  \end{tikzpicture}
  \caption{Problem formulation, problem solution, and interfaced solvers.}
  \label{fig:formulation-solution-solvers}
\end{figure*}

We turn to numerical results for power systems of various sizes.
We examine several combinations of the two problem formulations---feasibility~\eqref{eq:dist-feasibility-problem} and least-squares~\eqref{eq:dist-least-squares-problem}---and the two solution methods---\gls{admm} and \gls{aladin}, paired with different ways to compute sensitivities and interfance different solvers, see \autoref{fig:formulation-solution-solvers}.
The section is devised top-down: we begin with qualitative comparisons of \gls{admm} and \gls{aladin}, then examine the least-squares problem in combination with \gls{aladin} (for different solvers and different Hessian approximations).
The final section analyzes the convergence behavior for a 4662-bus system.

Our main finding is that the least-squares formulation with \gls{aladin} and a Gauss-Newton Hessian approximation outperforms all other combinations.

\begin{remark}[Settings common to all examples]
  \label{rema:settings}
  For all following examples, the connecting lines between all regions are modelled as transformers with a per-unit reactance of 0.00623, and a tap ratio of 0.985; the resistance, the total line charging susceptance, and the transformer phase shift angle are set to zero.\footnote{In light of \autoref{rema:wording} we switch back to referring to ``subsystems'' as ``regions'' and so on.}
  
  The initial condition for the primal state (i.e. the state of the electrical grid) is created from the \matpower case files as follows: the voltage angle and voltage magnitude are initialized with their respective entries from the entries in the \texttt{bus} struct; similaly, the net active power and the net reactive power are initialized as the difference between the respective summed entries in the \texttt{gen} struct and the \texttt{bus} struct.
  All dual variables are set to 0.01 initially.
  
  All computed solutions are verified relative to the reference solution provided by the \matpower command \texttt{runpf()}.
\end{remark}

\subsection{Qualitative comparison}
\label{sec:qualitative-comparison}
\begin{figure}
  \centering
    \begin{subfigure}{1\textwidth}
        \includegraphics[width=\textwidth]{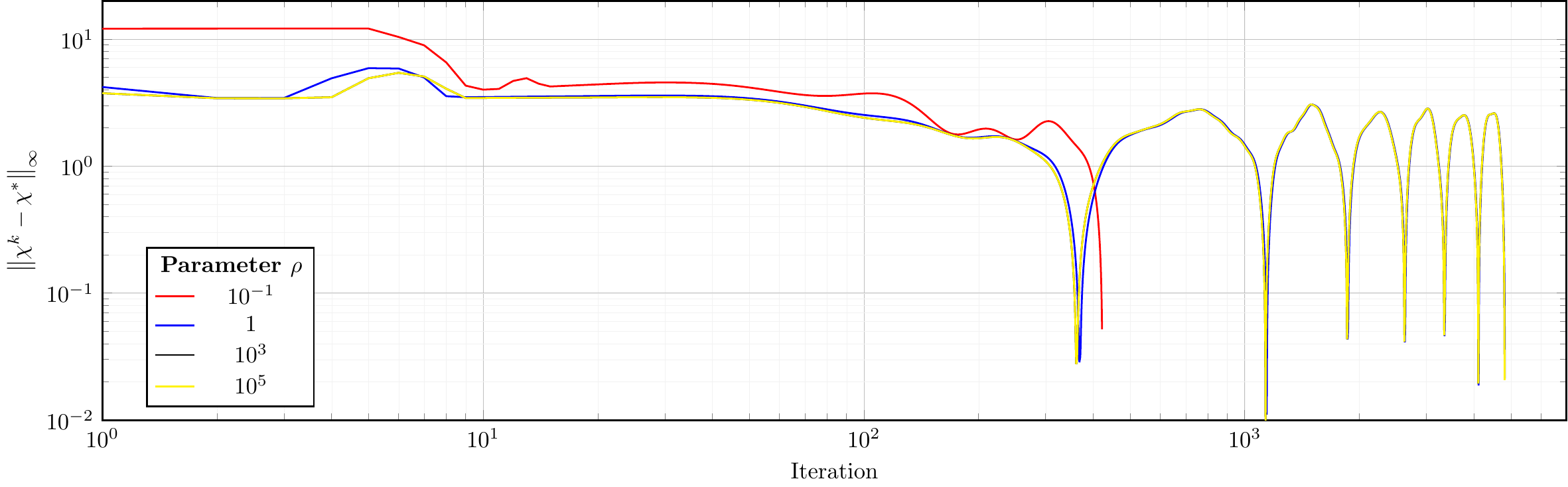}
        \includegraphics[width=\textwidth]{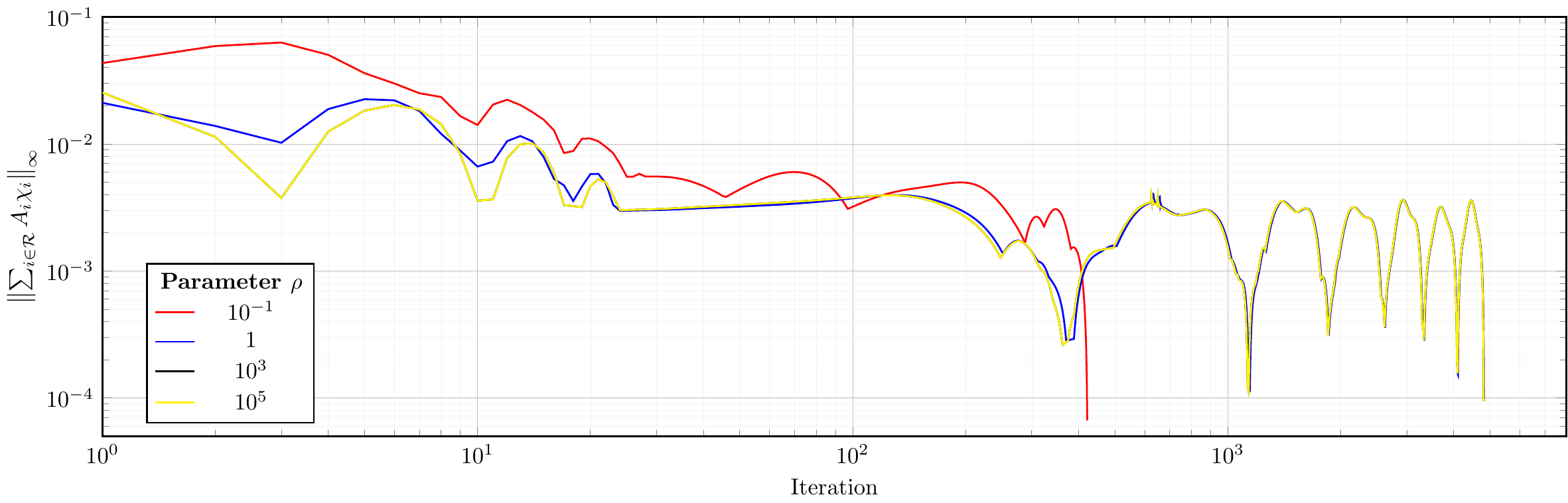}
        \caption{Varying penalty parameter~$\rho = 10^{n}$ for $n \in \{-1, 0, 3, 5 \}$.}
        \label{fig:admm-rho-dependence}
    \end{subfigure}
    \begin{subfigure}{1\textwidth}
      \includegraphics[width=\textwidth]{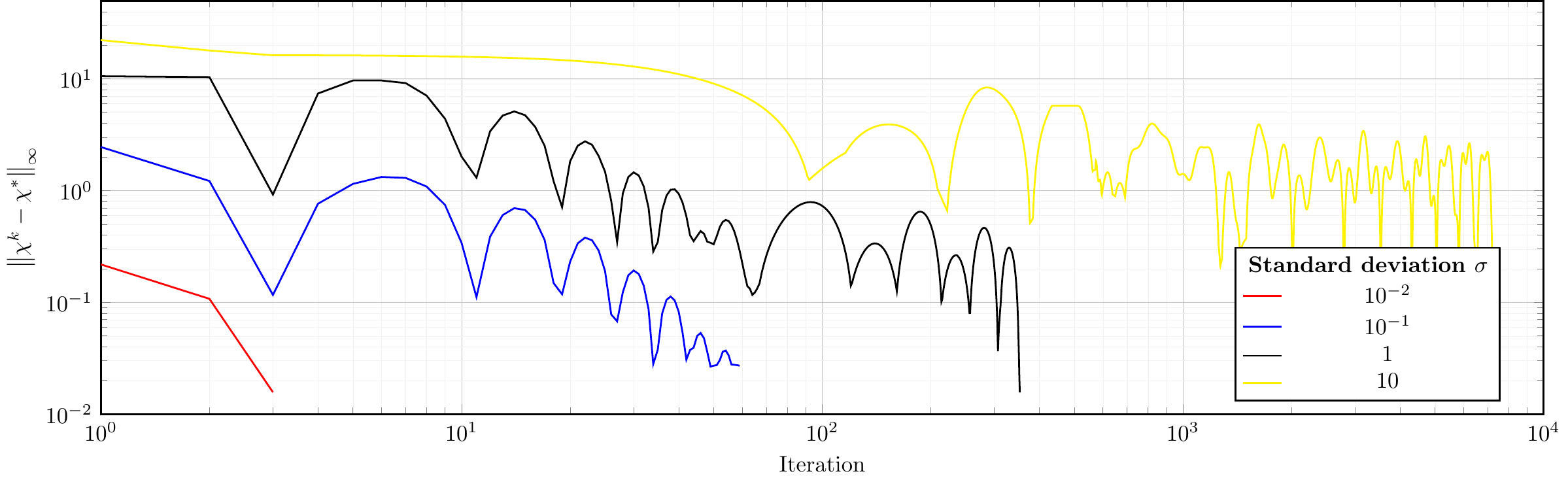}
      \includegraphics[width=\textwidth]{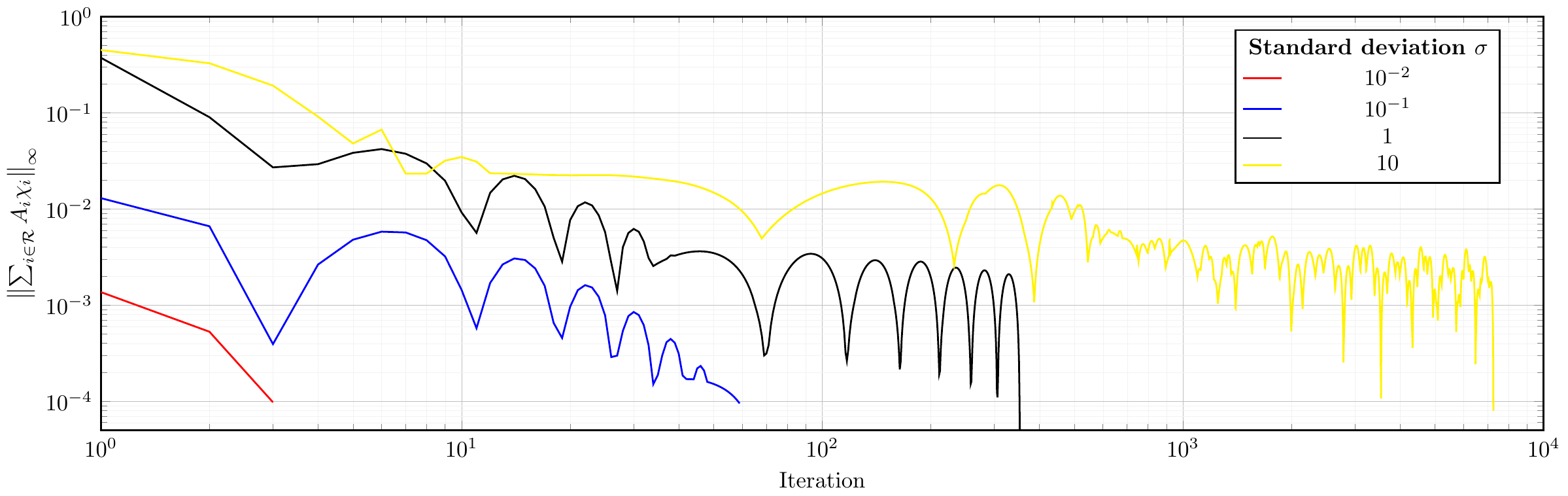}
      \caption{Varying intial conditions $\DistOptLocalState^\star + \alpha \hat{\DistOptLocalState}$, where $\DistOptLocalState^\star$ is the true solution, $\hat{\DistOptLocalState}$ is a vector whose entries are samples from a standard normal distribution, and $\sigma \in \{0.01, 0.1, 1, 10 \}$ is the standard deviation.}
      \label{fig:admm-initial-point-dependence}
    \end{subfigure}  
  \caption{Convergence behavior of \gls{admm} for a feasibility problem formulation of the 53-bus test case from \autoref{tab:computing-times}.
  In each subplot, the upper plot shows the distance to the optimal solution, and the lower plot shows the violation of the consensus constraints.}
\end{figure}


\begin{table}
  \centering
  \scriptsize
  \caption{Qualitative comparison of both problem formulations (feasibility vs. least-squares) and their solution by either \gls{admm} or \gls{aladin}.}
  \label{tab:qualitative-comparison}
  \renewcommand{\arraystretch}{1.3}
  \begin{tabular}{lcccg}
    \toprule
    & \multicolumn{2}{c}{Feasibility problem~\eqref{eq:dist-feasibility-problem}} & \multicolumn{2}{c}{Least-squares problem~\eqref{eq:dist-least-squares-problem}}\\
    & \gls{admm} & \gls{aladin} & \gls{admm}& \cellcolor{white}\gls{aladin} \\
    \midrule
    Scalability   & -{}-  & - & - & ++  \\
    Speed         & -{}-  & + & - & ++  \\
    Performance   & -{}-  & + & - & ++  \\
    Tuning        & -{}-  & - & - & +   \\
    \bottomrule
  \end{tabular}
\end{table}

We begin with a qualitative comparison of the applicability of both solution methods (\gls{admm} and \gls{aladin}) to both problem formulations (feasibility and least-squares); we base our qualitative findings on a total of 7 test cases that are summarized in the first four columns of \autoref{tab:computing-times}.

From \autoref{tab:qualitative-comparison}, which summarizes our qualitative findings, it appears that \gls{admm} is unsuitable for either problem formulation.
The performance of \gls{admm} depends critically on both the choice of the penalty parameter and the initial condition.
\autoref{fig:admm-rho-dependence} shows the convergence behavior for \gls{admm} applied to the feasibility formulation of the 53-bus test case from \autoref{tab:computing-times}.
For various choices of the penalty parameter~$\rho$ \gls{admm} exhibits strange and overall dissatisfying convergence properties.
Most of the considered cases (the ones from \autoref{tab:computing-times}) did not converge successfully even after having done significant parameter sweeps.
\autoref{fig:admm-rho-dependence} shows the influence of the choice of the penalty parameter~$\rho$, and the often-encountered convergence behavior with \gls{admm}: after relatively few iterations, the solution is in the vicinity of the optimal solution, but it takes several hundred iterations before further refinement occurs.
And even then, the solution is far from being accurate.
\autoref{fig:admm-initial-point-dependence} shows the critical dependence on the (primal) initial condition: Perturbing the primal initial condition around the opimal solution, the plots show that the entire optimization process is prolonged significantly.

In contrast to \gls{admm}, \gls{aladin} appears applicable to solve the distributed power flow problems from \autoref{tab:computing-times}.
In all qualitative aspects we consider (scalability, speed, performance, tuning), the least-squares formulation outperforms the feasibility counterpart by far, see \autoref{tab:qualitative-comparison}.
It is especially the aspect of scalability that hinders the feasibility problem: for instance, the 354-bus test case from \autoref{tab:computing-times} already took 38.2 s to solve with \fmincon, and converged within 14 \gls{aladin} iterations.

An explanation for this behavior may be the zero-cost objective function for the feasibility problem~\eqref{eq:dist-feasibility-problem}; sensitivities of the objective hence contain no information.
The advantage of the least-squares formulation is not just a non-zero objective function, but the absence of (in-)equality constraints in the local nonlinear programs, cf. \autoref{rema:local-optimization-problems} and \autoref{rema:nonlinear-least-squares}.

\subsection{Least-squares formulation with \gls{aladin}}

Based on our findings from the previous \autoref{sec:qualitative-comparison}, we consider only the least-squares formulation with \gls{aladin} in what is to follow.

\subsubsection{Different solvers}
\label{sec:different-solvers}

We investigate how the different solvers mentioned in \autoref{fig:rapidpf-flow-chart} cope with the different test cases from \autoref{tab:computing-times}; we use the sensitivities provided by \gls{rapidpf} in all cases, i.e. analytical gradients of the cost function, exact Jacobians, and the Gauss-Newton Hessian approximation.

Interestingly, \autoref{tab:computing-times} suggests that just half a dozen \gls{aladin} iterations are sufficient to solve the test cases, which range from a total of 53 buses to 4662 buses.
Hence, the applicability of \gls{aladin} itself is clearly demonstrated.
Of course, the overall solution time differs significantly with the choice of the local solver.\footnote{All computations were carried out on a standard a desktop computer with \textit{Intel\textsuperscript{\textregistered} Core\texttrademark\, i5-6600K CPU @ 3.50GHz} Processor and 16.0 \textsc{gb} installed \textsc{ram}; no efforts were made towards parallelization.\label{foot:computing-times}}
As a negative result we find that plain \aladinalpha, which interfaces only \casadi with Ipopt, is not suitable for the problem at hand.
That is why we chose to implement interfaces for the three other solvers: \fminunc, \fmincon, and \worhp.
Although \fminunc is the seemingly best fit---the local subproblems are unconstrained optimization problems---its practical applicability is limited to subproblems of a few hundred buses.
For the 2708- and 4662-bus test systems, \fminunc takes significantly longer, because the dimension of the local subproblem grows too large.
The solution times for \fmincon and \worhp are acceptable for all considered cases.
It stands to reason that \worhp will outperform \fmincon for even larger test cases, because it is able to exploit the sparsity of the optimization problem.

\begin{table}
  \centering
  \scriptsize
  \caption{Computing times for different test cases and different solvers when solving the distributed least-squares problem~\eqref{eq:dist-least-squares-problem} with \gls{aladin} and sensitivities from \gls{rapidpf}, see \autoref{foot:computing-times}.}
  \label{tab:computing-times}
  \rowcolors{2}{gray!10}{white}
  \begin{tabular}{rrp{1.2cm}rrrrrrr}
    \toprule
     &  & \matpower &  & \multicolumn{3}{c}{Solution time in s for} & \gls{aladin}\\
     \multirow{-2}{*}{Buses}& \multirow{-2}{*}{\nregions} & case files & \multirow{-2}{*}{\nconnections} & \fminunc & \fmincon & \worhp  & iterations\\
    \midrule
    53    &   3 & 9, 14, 30& 3 &   2.5 &  2.2  & 2.4 & 4\\ 
    354   &   3 & 3 $\times$ 118& 5&   2.5 & 3.1 & 4.8  & 5\\
    418   &   2 & 118, 300 & 2 &   4.5 & 5.2 & 7.0  & 5\\
    826   &   7 & 7 $\times$ 118 & 7 &   3.7 & 5.3 & 7.2  & 5\\
    1180  &   10 & 10 $\times$ 118 & 11 &   4.9 & 6.7 & 9.8  & 6\\
    2708  &   2 & 2 $\times$ 1354  & 1 &   212.7 & 41.9 & 53.6  & 4\\ 
    4662  &   5& 3 $\times$ 1354, 2 $\times$ 300& 4 &   387.9 &  90.1 & 113.8  & 5\\
    \bottomrule
  \end{tabular}
\end{table}

\subsubsection{Different Hessian approximations}
\label{sec:different-Hessian-approximations}
With \gls{aladin} being a second-order optimization method, the Hessian matrix is required---or an accurate yet easy-to-compute approximation thereof.
We compare four different Hessian approximations for the least-squares problem~\eqref{eq:dist-least-squares-problem} with \gls{aladin}: finite differences, \gls{bfgs}, limited-memory \gls{bfgs}, and the Gauss-Newton method.
The results, which are shown in \autoref{tab:Hessian-approximations}, confirm what is to be expected: the Gauss-Newton method outperforms all other methods.
This is in accordance with the fact that exploiting the structure of the least-squares formulation correctly pays off tremendously.
The finite difference approximation, just like the two \gls{bfgs} methods, are all-purpose Hessian approximation unaware of the underlying problem structure.
Gauss-Newton, in turn, is a Hessian approximation tailored to nonlinear least squares problem, see also \autoref{rema:sensitivities}.
The results from \autoref{tab:Hessian-approximations} make it clear that already for small system sizes, the all-purposes Hessian approximations should be avoided, because they lead to longer computation times.\footnote{Note however that the total number of \gls{aladin} iterations is unaffected.}
Consequently, the default Hessian approximation for least-squares problems is the Gauss-Newton method in \gls{rapidpf}.

\begin{table}[h!]
  \centering
  \scriptsize
  \caption{Computing times for least-squares problem~\eqref{eq:dist-least-squares-problem} with \gls{aladin} and \fmincon, for different Hessian approximations.
  The entries in the column ``Buses'' refers to the entries in \autoref{tab:computing-times}.}
  \label{tab:Hessian-approximations}
  \rowcolors{2}{gray!10}{white}
  \begin{tabular}{rrrrr}
    \toprule
    Buses & Finite difference & \gls{bfgs} & Limited-memory \gls{bfgs} & Gauss-Newton \\
    \midrule
    53    &   10.0  &   28.6    &   22.9  &   2.2 \\
    354   &   61.5  &   287.8   &   107.4 &   3.1 \\
    418   &   185.6 &   1086.4  &   148.2 &   5.2 \\
    826   &   n/a   &   n/a     &   n/a   &   5.3 \\
    $\hdots$ & n/a   &   n/a     &   n/a   &   See \autoref{tab:computing-times} \\
    4662 & n/a   &   n/a     &   n/a   &   See \autoref{tab:computing-times} \\
    \bottomrule
  \end{tabular}
\end{table}

\subsection{4662-Bus system -- Convergence behavior}
Next, we study the convergence behavior of the 4662-bus test case.
This test case is composed of three 1354-bus \matpower test cases, and two 300-bus \matpower test cases.
\autoref{tab:4662-bus-test-case} shows the connecting buses between the regions.
For other relevant information such as how the connecting lines are modelled, and how the initial conditions are chosen, see \autoref{rema:settings}.

\begin{figure}
  \centering
    \begin{subfigure}{1\textwidth}
        \includegraphics[width=\textwidth]{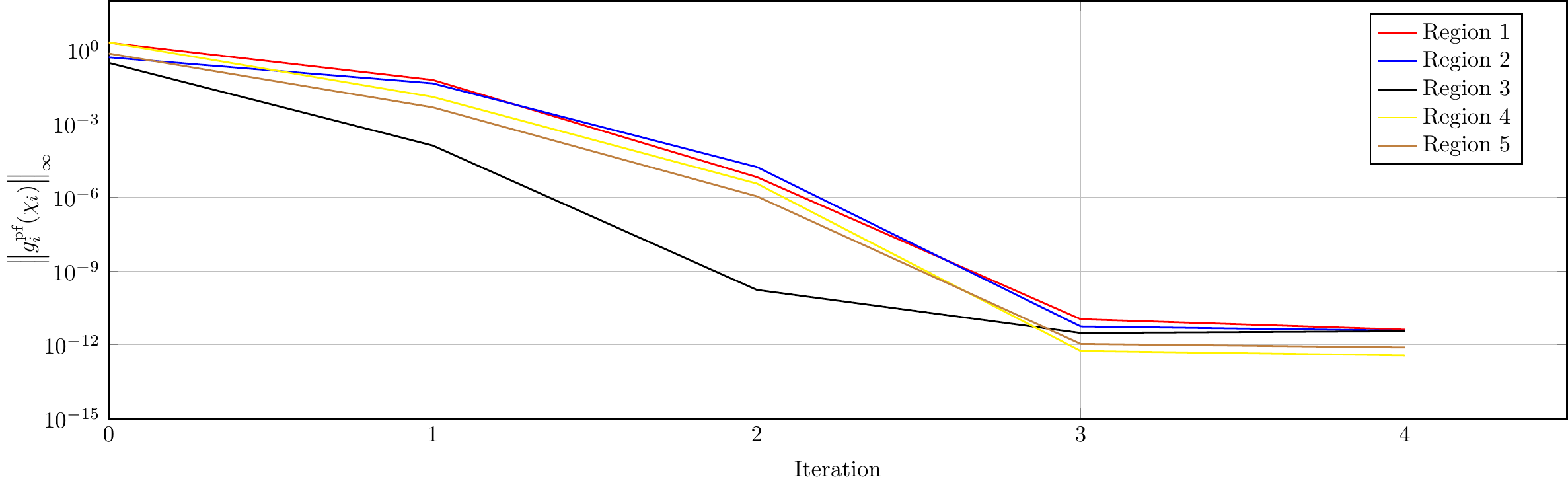}
    \end{subfigure}
    \begin{subfigure}{1\textwidth}
        \includegraphics[width=\textwidth]{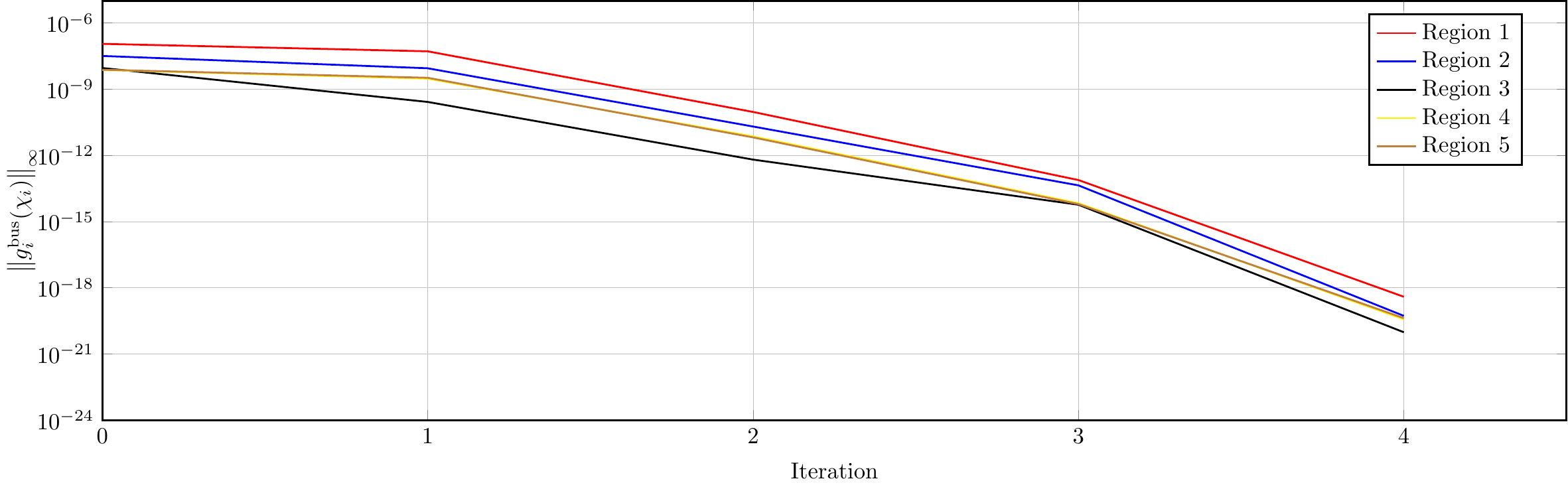}
    \end{subfigure} 
    \begin{subfigure}{1\textwidth}
        \includegraphics[width=\textwidth]{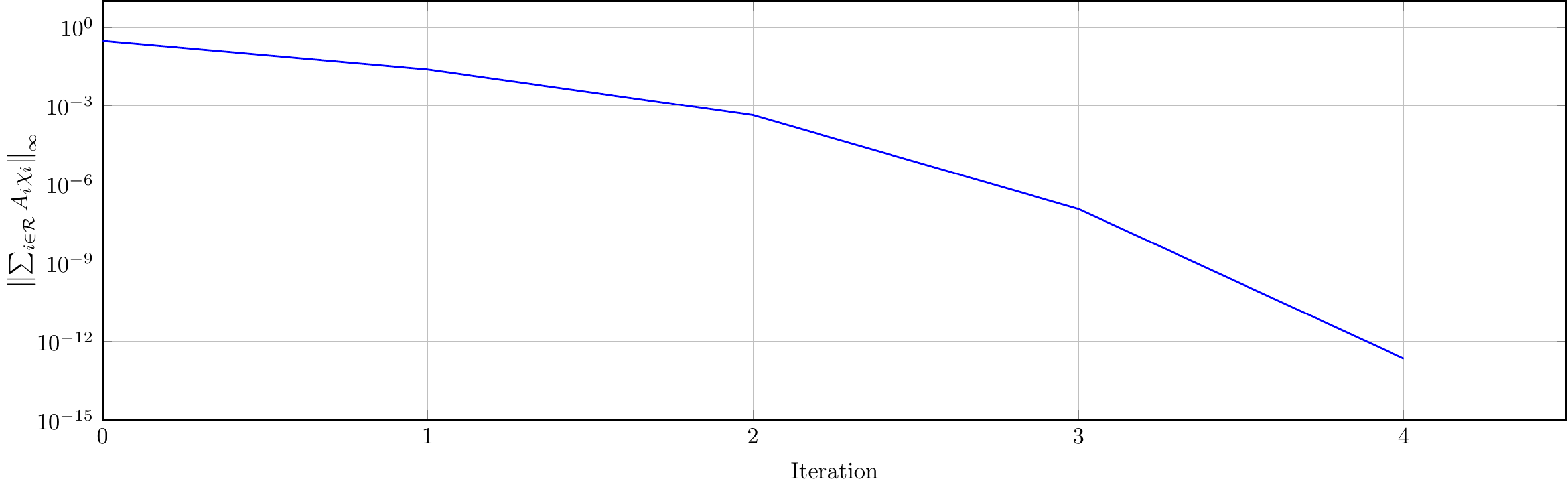}
    \end{subfigure}  
  \caption{Decrease of the $\infty$-norm of the power flow equations, the bus specifications, and the consensus violations, each per \gls{aladin} iteration for the 4662-bus system from \autoref{tab:computing-times}.}
  \label{fig:4662-violations}
\end{figure}

To solve the distributed power flow problem we choose a least-squares formulation with \gls{aladin}.
We use the Gauss-Newton Hessian approximation, and \fmincon is used to solve the local problems.
From \autoref{tab:computing-times} we see that this setup requires 5 \gls{aladin} iterations and about 90 seconds.
\autoref{fig:4662-violations} shows, for every \gls{aladin} iteration and for every region, the $\infty$-norm of the power flow equations~\eqref{eq:dist-power-flow-problem-pf}, of the bus specifications~\eqref{eq:dist-power-flow-problem-busspecs}, and of the consensus constraint violations~\eqref{eq:dist-power-flow-problem-consensus}.
After 5 \gls{aladin} iterations, all violations are below $10^{-10}$, and the computations are terminated.

\begin{table}
  \centering
  \scriptsize
  \caption{Regions and used test cases for 4662-bus test case (left). Connecting buses between regions (middle). Connection graph (right)}
  \label{tab:4662-bus-test-case}
  \begin{tabular}{ccc}
    \begin{tabular}{rr}
      \toprule
      Region & \matpower case file \\
      \midrule
      1 & \texttt{case1354pegase} \\
      2 & \texttt{case1354pegase} \\
      3 & \texttt{case1354pegase} \\
      4 & \texttt{case300} \\
      5 & \texttt{case300} \\
      \bottomrule
    \end{tabular}
    &
    \begin{tabular}{rrrr}
      \toprule
      \multicolumn{2}{c}{From-system} & \multicolumn{2}{c}{To-system} \\
      Region & Bus & Region & Bus \\
      \midrule
      1 & 17 & 2 & 46 \\
      1 & 111 & 3 & 271 \\
      2 & 64 & 4 & 10 \\
      2 & 837 & 5 & 8 \\
      \bottomrule
    \end{tabular}
    &
    \begin{tabular}{c}
      \begin{tikzpicture}
        [ grow = down,
          edge from parent/.style = {draw, -latex},
          every node/.style = {font=\scriptsize, circle, draw},
          sibling distance = 4em,
          level distance = 4em,
          inner sep = 2pt,
        ]
        \node {$\mathcal{R}_1$}
          child { node {$\mathcal{R}_3$} }
          child { node {$\mathcal{R}_2$}
            child { node {$\mathcal{R}_4$} }
            child { node {$\mathcal{R}_5$}}
             };
      \end{tikzpicture}
    \end{tabular}
  \end{tabular}
\end{table}

\section{Conclusion \& Outlook}
\label{sec:conclusion}
The relevance of distributed power flow problems is increasing, because their solutions allow for better cooperation between different stakeholders, e.g. \gls{tsos} and \gls{dsos}.
Distributed optimization is a viable technique to tackle such distributed power flow problems.
It is speficically the \acrfull{aladin} with its convergence guarantees that yields promising results: if the distributed power flow problem is formulated as a distributed least-squares problem, and if a Gauss-Newton Hessian approximation is used, then about half a dozen iterations suffice to converge to the correct solution.
To facilitate rapid prototyping we introduce \acrfull{rapidpf}, which is fully \matpower-compatible \matlab code that takes over the laborious task of creating code amenable to distributed optimization.

Future steps will focus mainly on further structure exploitation for solving the problem, and on implementing larger test cases.
The least-squares formulation is promising, hence further improvements are possible, such as relying not on an all-purpose solver but devising a solver dedicated to nonlinear-least squares problems.
A first step might be a tailored Gauss-Newton method, or a tailored Levenberg-Marquardt method~\cite{Nocedal2006}.
For the Gauss-Newton method, for example, it is possible to avoid having to compute the Hessian altogether, because a singular-value decomposition or a conjugate gradient method can be applied directly to solve the linearized problem~\cite{Nocedal2006}.
This user-defined nonlinear-least squares solver must then be interfaced with \aladinalpha.

The simulation results we presented are all carried out on a single machine.
To leverage the literal \emph{distribution} of the optimization, efforts toward parallel computing shall be undertaken when tackling larger test cases.

Finally, \gls{rapidpf} can be extended to optimal power flow problems upon adding cost functions per region.

\section*{Acknowledgment}
The authors would like to thank Jochen Bammert and Tobias Wei\ss bach (both TransnetBW GmbH) for insightful discussions and continuing support and interest in distributed power flow.
Finally, Tillmann M\"uhlpfordt thanks Daniel Bacher for his supporting the migration from Gitlab to GitHub.

\section*{Contributor roles}
See \autoref{tab:contributor-roles} for the roles of each author.
\begin{table}[h!]
  \scriptsize
  \centering
  \caption{Contributor roles.}
  \label{tab:contributor-roles}
  \begin{tabular}{lp{0.57\textwidth}}
    \toprule
    Author & Role(s) \\
    \midrule
    Tillmann M\"uhlpfordt & Conceptualization, investigation, methodology, software, writing -- original draft\\
    Xinliang Dai & Software, investigation, visualization\\
    Alexander Engelmann & Methodology (while at \textsc{kit}), writing -- original draft (\autoref{sec:solution-admm}, \autoref{sec:solution-aladin}) (while at \textsc{tu} Dortmund)\\
    Veit Hagenmeyer & Conceptualization, funding acquisition, writing -- review \& editing\\
    \bottomrule
  \end{tabular}
\end{table}

\footnotesize
\bibliography{lit, lit_old, alex}

\end{document}